\documentclass[]{article} 

\usepackage{graphicx}
\usepackage{array}
\usepackage{amsfonts}
\def\norm#1{\left\Vert#1\right\Vert}

\def\eps{\varepsilon}                                                         
\def\crm{\cr\noalign{\medskip}}
\def\abs#1{\left\vert#1\right\vert}
\def\pabs#1{\vert#1\vert}
\def\m@th{\mathsurround=0pt}
\def\EQM#1{\vcenter{\normalbaselines\m@th
    \ialign{${\displaystyle ##}$\hfil&&\ ${\displaystyle ##}$\hfil\crcr
    \mathstrut\crcr\noalign{\kern-\baselineskip}
    \noalign{\smallskip}
    #1\crcr\mathstrut\crcr\noalign{\kern-\baselineskip}}}}

\def\om{\omega}
\def\Om{\Omega}
\def\s{\sigma}

\def\vphi{\varphi}
\def\s{\sigma}
\def\r{\right}
\def\l{\left}

\def\tx{{\tilde x}}
\def\e{{\rm e}} 
\def\trait{\noalign{\smallskip\hrule\smallskip}}
\def\bbbn{\mathbb{N}}
\def\bbbr{\mathbb{R}}
\def\bbbt{\mathbb{T}}
\def\bbbc{\mathbb{C}}
\def\bbbz{\mathbb{Z}}

\def\llabel{\label} 
\newcommand{\Frac}[2]{{{\displaystyle\strut#1}\over{\displaystyle\strut#2}}}
\def\eps{\varepsilon}
\def\etal{{\it et al.}} 
\def\abs#1{\left\vert#1\right\vert}
\def\om{\omega}
\newcommand{\FFrac}[2]{{{\displaystyle\strut#1}\over{\displaystyle\strut#2}}}

\newcommand{\Dron}[2]{\FFrac{\partial#1}{\partial#2}}
\newcommand\be{\begin{equation}}
\newcommand\ee{\end{equation}}
\def\bib#1#2#3#4#5#6#7{\bibitem{#1} {#2}\ :{\ #3},\ {\ #4},\ {{\it #5}},\
{{\bf #6}},\ {(#7)} \par }
\def\bibC#1#2#3#4#5#6{\bibitem{#1} {#2}\ :{\ #3},\ {\ #4},\ {{\it #5}},\
{(#6)} \par }

\newcommand\figpath{}
\newcommand\figa{\figpath 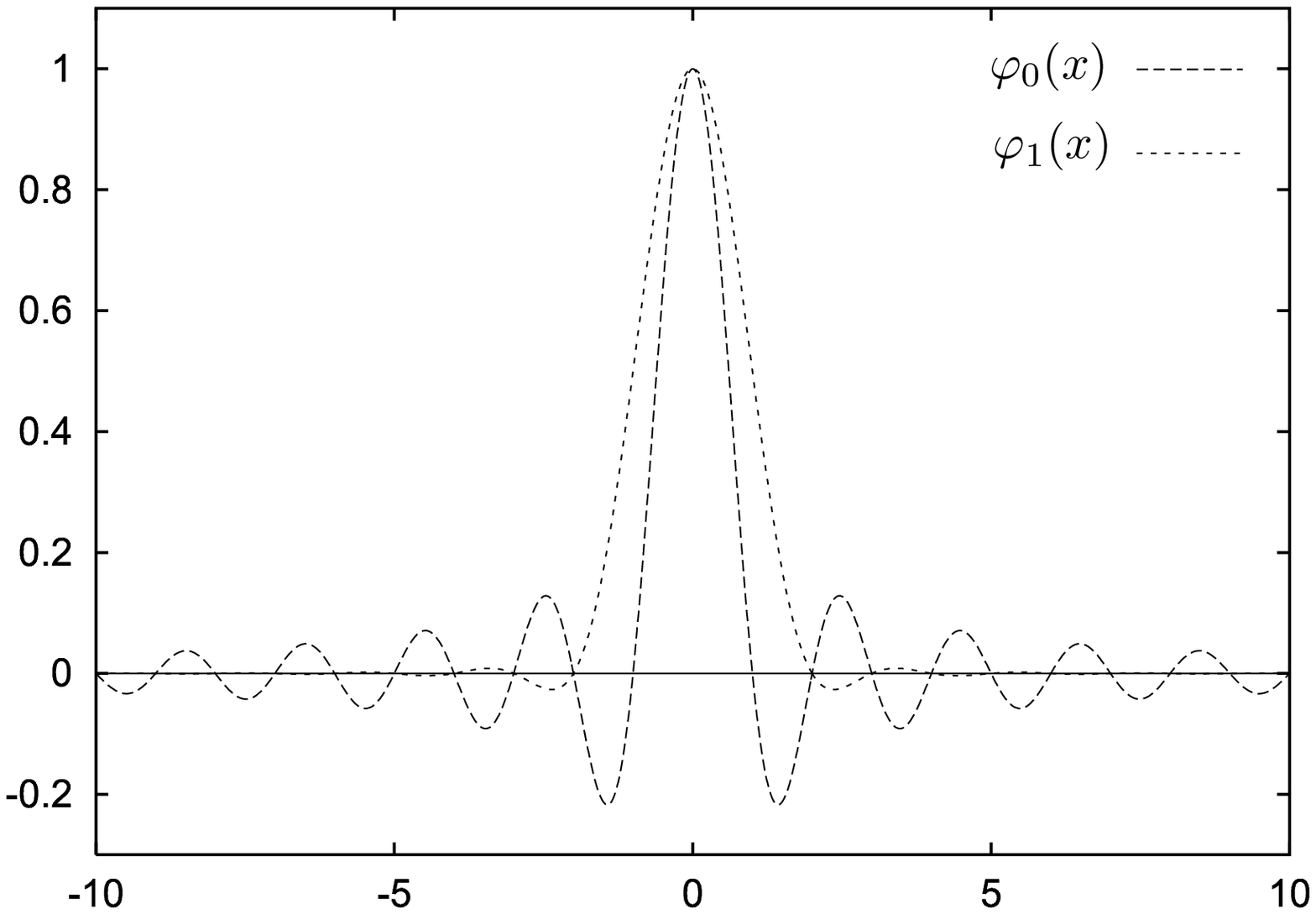}
\newcommand\figb{\figpath 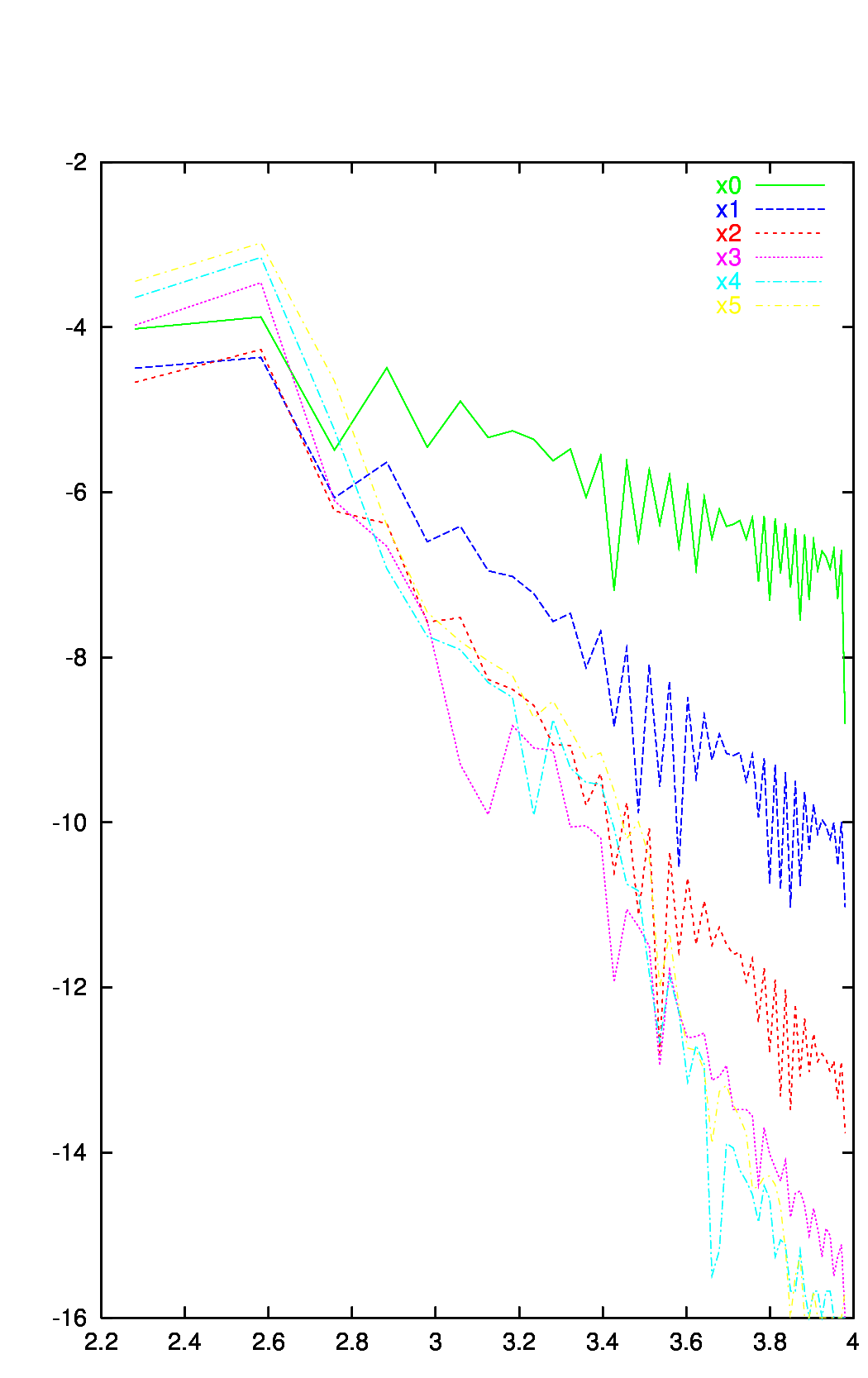}
\newcommand\figc{\figpath 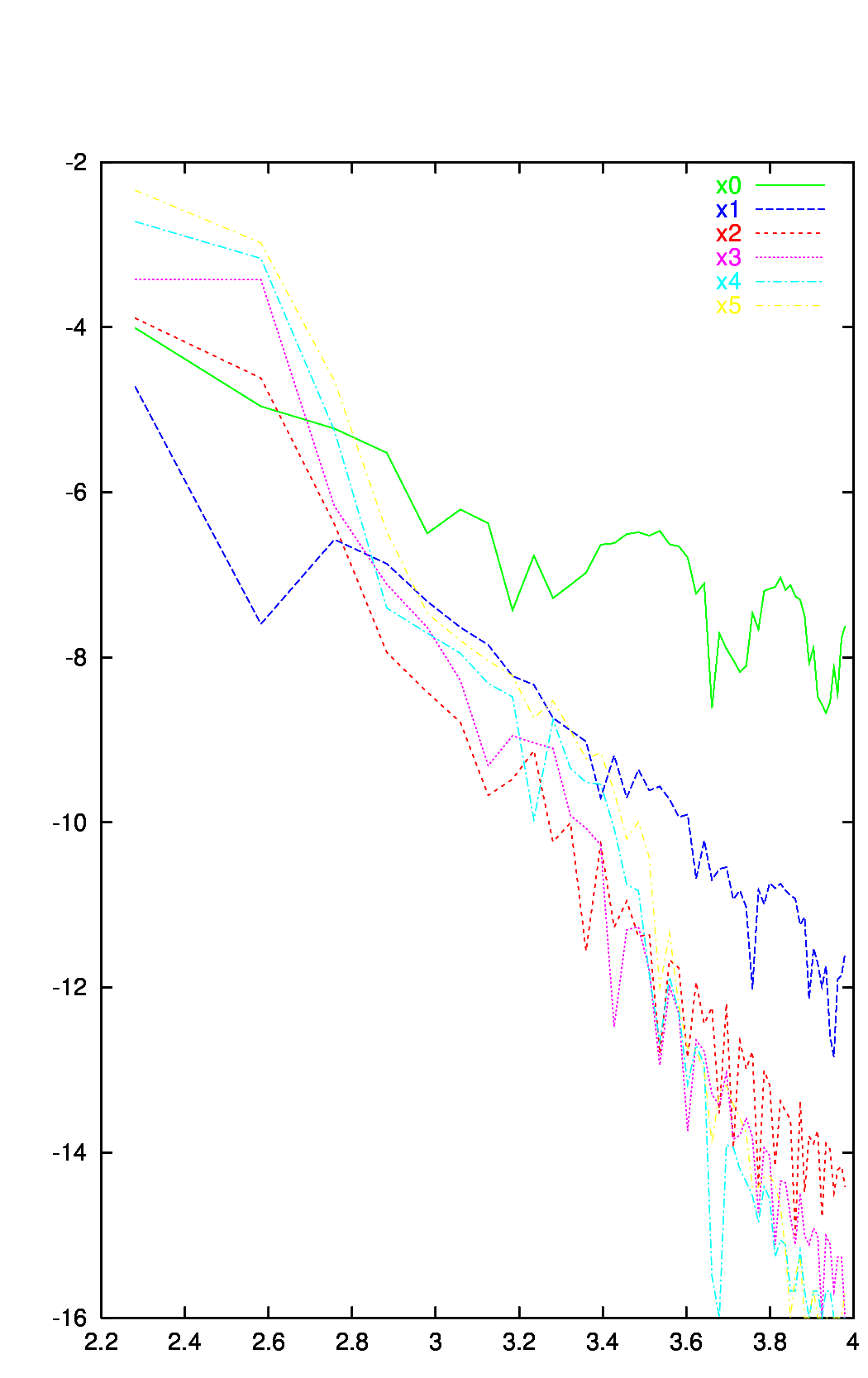}
\newcommand\figd{\figpath 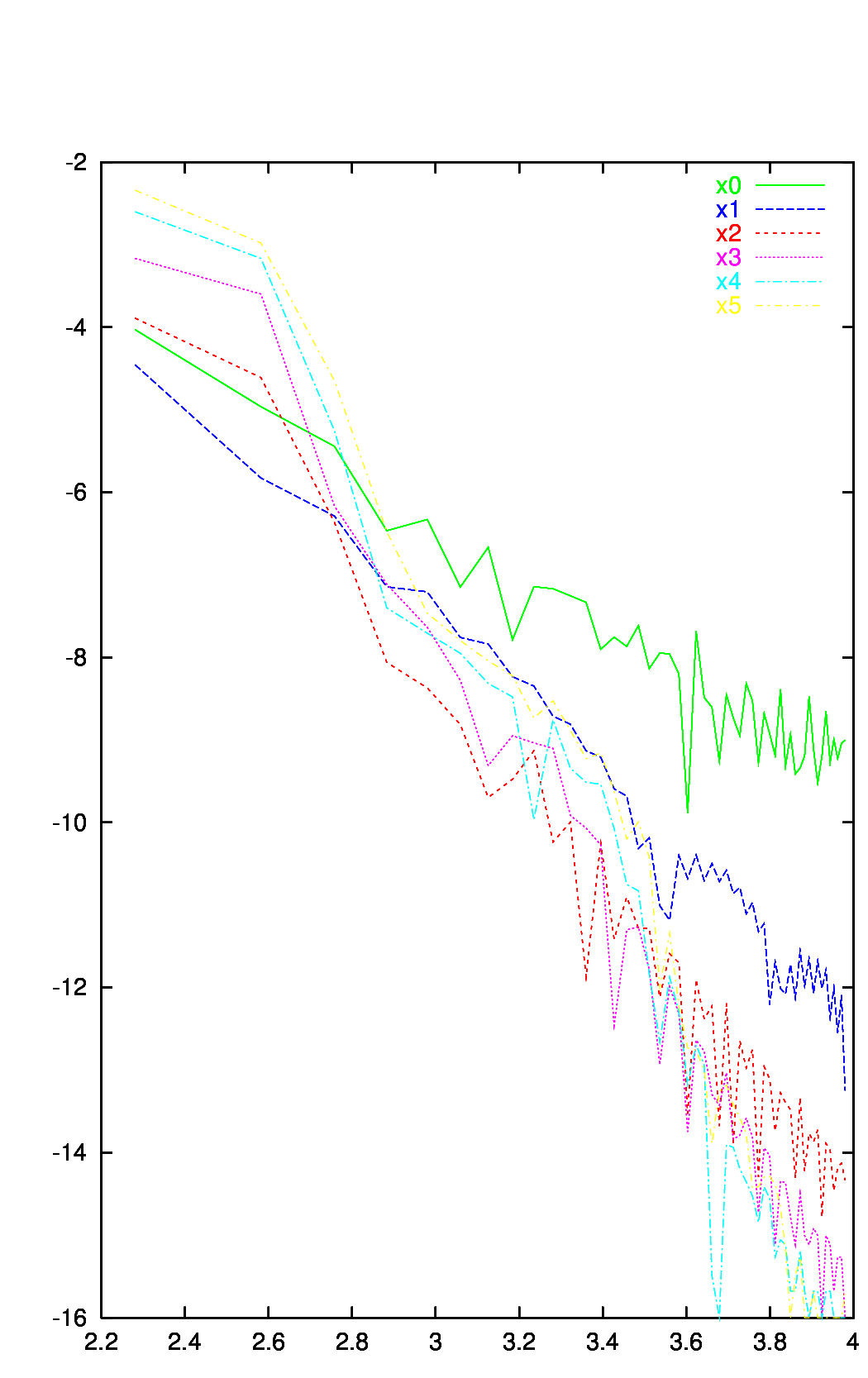}
\newcommand\fige{\figpath 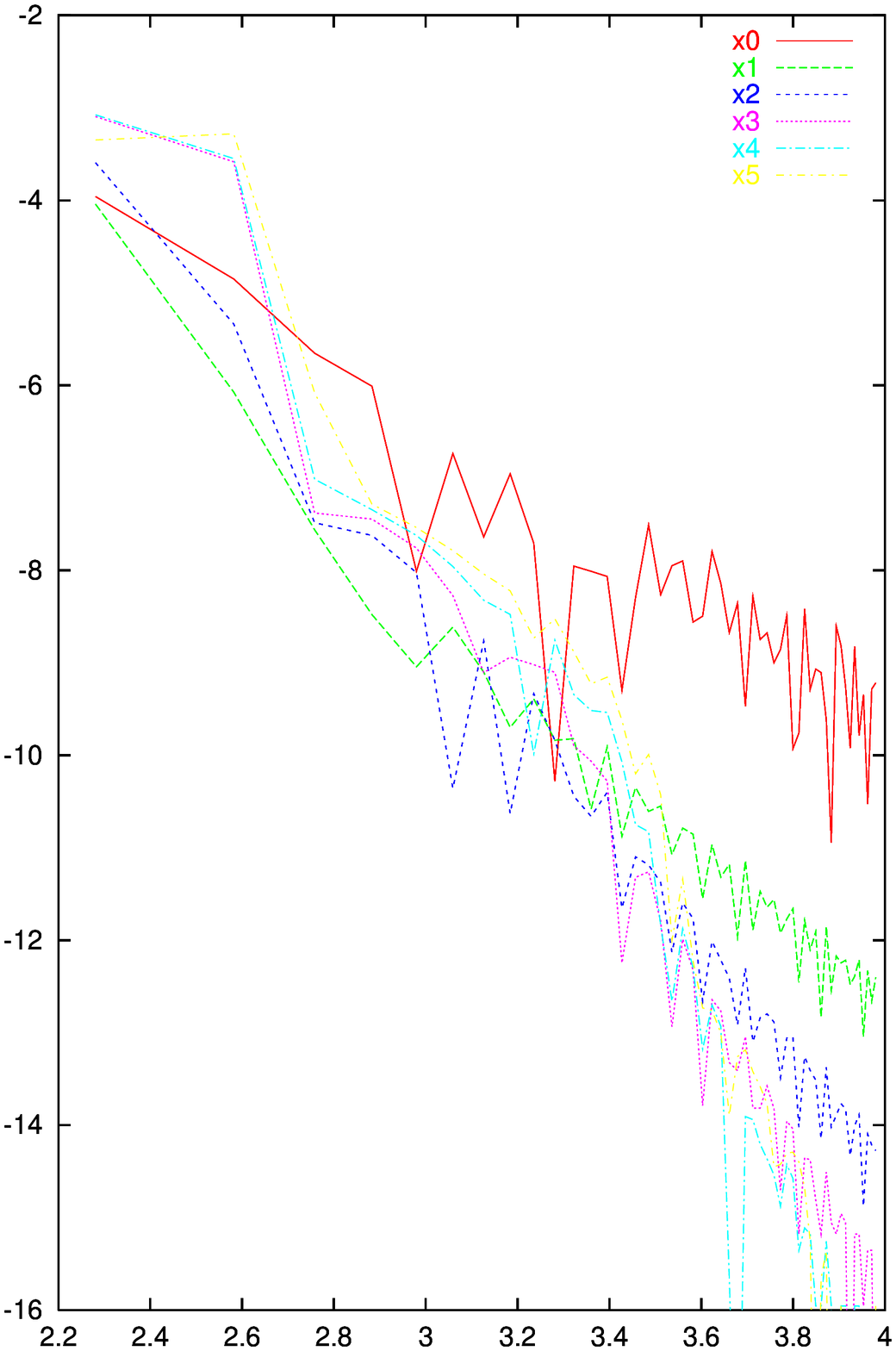}
\newcommand\figf{\figpath 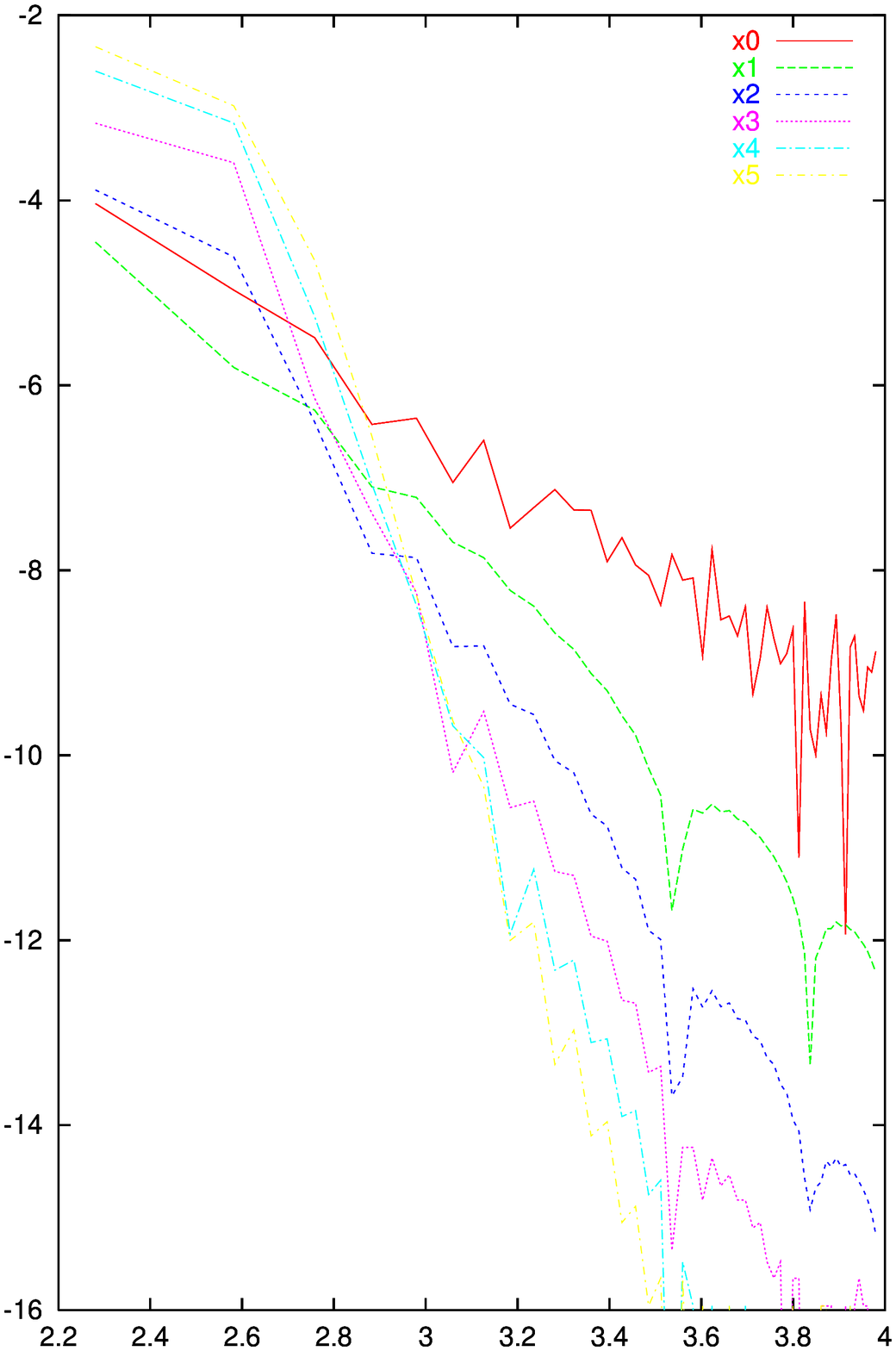} 
\newcommand\figg{\figpath 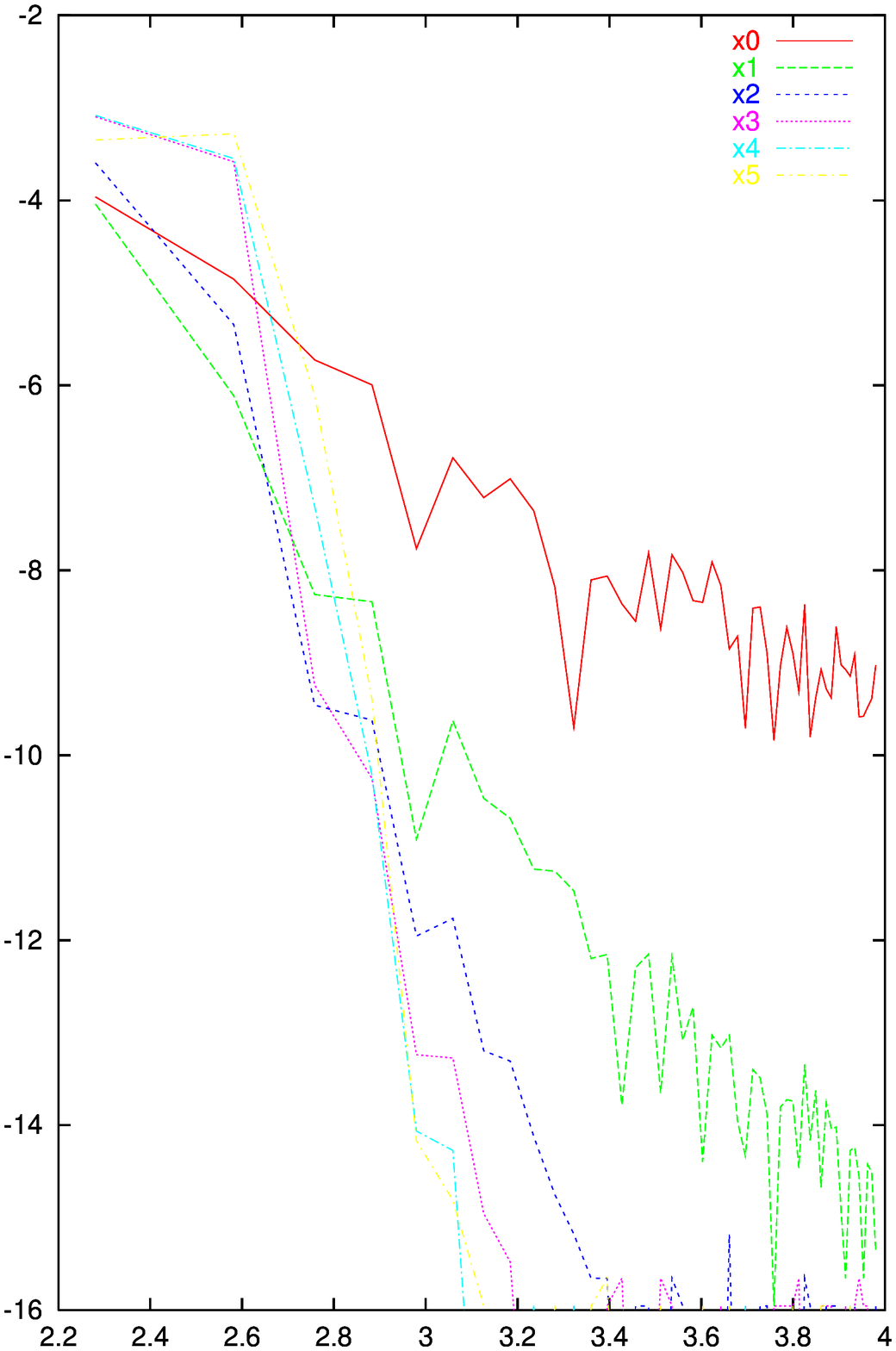}

\newtheorem{lemme}{Lemma.}
\newtheorem{definition}{Definiton.}
\newtheorem{cor}{Corollary.}
\newtheorem{theo}{Theorem.}
\newtheorem{prop}{Proposition.}

\title{Frequency map analysis and quasiperiodic decompositions.} \author{Jacques
Laskar\\ \it \small Astronomie et Syst\`emes Dynamiques, \\ 
\it \small CNRS UMR8028, IMCCE-Observatoire de Paris,\\
 \it \small 77 Av. Denfert-Rochereau, 75014 Paris, France}

\date{13 march 2003, revised 26 june 2003 \\ proceedings of Porquerolles School, sept. 2001}

\begin{document}

\maketitle

\section{Introduction}
Frequency Map Analysis is a numerical method based on refined Fourier techniques 
which provides a clear representation of the global dynamics of many multi-dimensional 
systems, and which is particularly adapted for systems of  3-degrees of freedom and more.
This method relies heavily on the possibility of making accurate 
quasiperiodic approximations of of  quasiperiodic signal 
given in a numerical way. 
In the present paper, we will describe the basis of the frequency analysis method, 
focussing on the quasi periodic approximation techniques.
Application of these methods for the study of the global dynamics 
and chaotic diffusion of Hamiltonian systems and symplectic maps  in different domains 
can be found in (Laskar, 1988, 1990, Laskar and Robutel, 1993, Robutel and Laskar, 2001,
Robutel, 2003, this volume, and
Nesvorn\'y and Ferraz-Mello, 1997) for solar system dynamics, and  in 
(Papaphilippou and Laskar, 1996, 1998, Laskar, 2000, Wachlin and Ferraz-Mello,Ê1998, 
Valluri and Merritt, 1998, Merritt and Valluri, 1999) for
galactic dynamics. The method has been particularly successful for its application in particle 
accelerators (Dumas and Laskar, 1993, Laskar and Robin, 1996, Robin \etal, 2000, Comunian \etal,
2001,  Papaphilippou and Zimmermann, 2002, Steier \etal, 2002), and   was also  used for the understanding of atomic
physics (Milczewski \etal, 1997),  or more general dynamical system issues (Laskar \etal, 1992,
Laskar, 1993, 1999,  Chandre \etal, 2001).

\section{Frequency Maps}
According to the KAM theorem (see Arnold \etal, 1988), in the
phase space of a sufficiently close to integrable conservative system, many
invariant tori will persist. Trajectories starting on one of these tori remain
on it thereafter, executing quasiperiodic motion with a fixed frequency vector
depending only on the torus. The family of tori is parameterized over a Cantor
set of frequency vectors, while in the gaps of the Cantor set chaotic behavior
can occur.
The frequency analysis algorithm  will
numerically compute over a finite time span a frequency vector for any initial
condition. On the KAM tori, this frequency vector will be a very accurate
approximation of the actual frequencies, while in the chaotic regions, 
the algorithm will still provide some determination of the frequency vector, 
but  in this region, complementary of the KAM tori,
the frequency vector  will not be uniquely defined.

Let us consider a $n$--DOF Hamiltonian system close to integrable in the form $
H( I,\theta) = H_0(I) + \eps H_1( I,\theta) \label{hiteta} $, where $H$ is real
analytic for $( I,\theta)\in B^n\times \bbbt^n$, $B^n$ is a domain of $\bbbr^n$
and $\bbbt^n$ is the $n$-dimensional torus. For $\eps = 0$, the Hamiltonian
reduces to $H_0(I) $ and is integrable. The equations of motion are then
for all $ j=1,\dots,n$
\begin{equation}
\dot I_j = 0\;\;, \quad \dot
\theta_j = \Dron{H_0(I)}{I_j}=\nu_j(I)\ ;
\end{equation}
which gives $z_j(t) =z_{j0} e^{i\nu_j t} \label{zj} $ in the complex variables
$z_j = I_j\exp i\,\theta_j$, where $z_{j0} = z_j(0)$. The motion in phase
space takes place on tori, products of true circles with radii
$I_j=\abs{z_j(0)}$, which are described at constant velocity $\nu_j(I)$. If
the system is nondegenerate, that is if
\begin{equation}
det\pmatrix{\Dron{\nu(I)}{I}\cr}=det\pmatrix{\Dron{{}^2H_0(I)}{I^2}\cr}\neq 0
\label{eq.torsion}
\end{equation}
the frequency map $F: B^n \longrightarrow \bbbr^n; \quad (I) \longrightarrow
(\nu) \quad$ is a diffeomorphism on its image $\Omega$, and the tori are as
well described by the action variables $(I) \in B^n$ or in an equivalent manner
by the frequency vector $(\nu)\in \Omega$.  For a nondegenerate system,
KAM theorem
still asserts that for sufficiently small values of $\eps$ , there exists a
Cantor set $\Omega_\eps$ of values of $(\nu)$, satisfying a Diophantine
condition of the form
\be
\abs{(k,\nu)} >{\kappa_\eps}/{\abs{k}^m}
\label{eq.dio}
\label{eq.3}
\ee
for which the perturbed system still possesses smooth invariant tori with
linear
flow (the KAM tori). Moreover, according to P\"oschel (1982), there exists a
diffeomorphism
\begin{equation}
\Psi: \bbbt^n \times \Omega \longrightarrow \bbbt^n\times B^n; \qquad
(\varphi, \nu)
\longrightarrow (\theta,I)
\end{equation}
which is analytical with respect to $\varphi$, $C^\infty$ in $\nu$, and on
$\bbbt^n\times \Omega_\eps$ transforms the Hamiltonian equations into the trivial system
\be
\dot \nu_j = 0\;, \quad \dot \varphi_j = \nu_j \ .
\ee
 For frequency vectors
$(\nu)$ in $\Omega_\eps$, the solution lies on a torus and is given in complex
form by its Fourier series
\begin{equation}
z_j(t) =z_{j0} e^{i\nu_j t} + \sum_{m} a_{m}(\nu) e^{i<m,\nu>t}
\label{zetan}
\label{eq.5}
\end{equation}
where the coefficients $a_{m}(\nu)$ depend smoothly on the frequencies $(\nu)$.
If we fix $\theta \in \bbbt^n$  to some value $\theta=\theta_0$, we obtain a
frequency map on $B^n$ defined as
\begin{equation}
 F_{\theta_0}: B^n\longrightarrow \Omega\;\;; \qquad I \longrightarrow
p_2(\Psi^{-1}(\theta_0,I))
\label{fremap}
\end{equation}
where $p_2$ is the projection on $\Omega$ ($p_2(\phi,\nu)=\nu$).
For sufficiently small $\eps$, the torsion condition
(\ref{eq.torsion}) ensures that the frequency map $F_{\theta_0}$ is a smooth
diffeomorphism.

\subsection{Isoenergetic nondegeneracy}

If the isoenergetic nondegeneracy condition of Arnold (1989) is verified
\be
det\pmatrix{\Dron{\nu(I)}{I} & \nu(I)\cr \nu(I)^\tau & 0 \cr}=
det\pmatrix{\Dron{{}^2H_0(I)}{I^2} & \Dron{H_0(I)}{I}\cr \Dron{H_0(I)}{I}^\tau & 0}\neq 0
\ee
then, if $\nu_n=\Dron{H_0(I)}{I_n}\neq 0$,  and for small enough $\eps$, the application 
\be
(I_1,\dots,I_{n-1},I_n)\longrightarrow \left(\Frac{\nu_1}{\nu_n},\dots,\Frac{\nu_{n-1}}{\nu_n},H\right)
\ee
is a diffeomorphism, and so will be its restriction $F'_{\theta_0}$ on $H^{-1}(h)$
\be
(I_1,\dots,I_{n-1})\longrightarrow \left(\Frac{\nu_1}{\nu_n},\dots,\Frac{\nu_{n-1}}{\nu_n}\right) 
\ee
When the isoenergetic nondegeneracy condition is verified, 
 we can thus restrict ourself to an energy level $H=h$, 
 and  the frequency map $F_{\theta_0}:\bbbr^n \rightarrow \bbbr^n$ 
can be reduced to a map $F'_{\theta_0}:\bbbr^{n-1} \rightarrow \bbbr^{n-1}$.
\section{Quasiperiodic approximations}

The frequency analysis method and algorithms rely heavily on the observation that when a
quasiperiodic function $f(t)$ in the complex domain $\bbbc$ is given
numerically, it is possible to recover a quasiperiodic approximation of $f(t)$
in a very precise way over a finite time span $[-T,T]$, several orders of
magnitude more precisely than  by simple Fourier analysis.  Indeed,
let
\begin{equation}
f(t) = e^{i \nu_1 t} + \sum_{k\in \bbbz^n-(1,0,\dots,0)} a_k e^{i \langle
k,\nu \rangle t}\;\; ,
\quad a_k\in
\bbbc
\label{eq.ft}
\end{equation}
be a KAM quasiperiodic solution of an Hamiltonian system in $B^n\times
\bbbt^n$, where the frequency vector $(\nu)$ satisfies a Diophantine condition
(\ref{eq.dio}). The frequency analysis algorithm NAFF will provide an
approximation $f'(t) = \sum^{N}_{k=1} a'_k e^{i\om'_k t}$  of $f(t)$ from its
numerical knowledge over a finite time span $[-T,T]$ . The frequencies $\om_k'$
and complex amplitudes $a_k'$ are computed through an iterative scheme. In
order to determine the first frequency $\om_1'$, one searches for the maximum
amplitude of $\phi(\sigma) = \langle f(t), e^{i\sigma t}\rangle$ where the
scalar product $\langle f(t), g(t)\rangle$ is defined by
\begin{equation}
\langle f(t), g(t)\rangle = \Frac{1}{2T}\int_{-T}^T f(t)\bar g(t) \chi(t)
dt \ ,
\end{equation}
and where $\chi(t)$ is a weight function (see next section). Once the first
periodic term $e^{i\om_1' t}$ is found, its complex amplitude $a_1'$ is
obtained by orthogonal projection, and the process is restarted on the
remaining part of the function $f_1(t) = f(t) -a_1' e^{i\om_1' t}$. 

In the next sections, we provide the rigorous foundations 
of the frequency analysis algorithm by  showing that this algorithm converges 
towards  the frequency map described formally in the previous sections. 
It is interesting  to see that the convergence of the 
frequency map algorithm will require some conditions 
on the frequency vector that will always be satisfied in 
the case of a KAM regular solution.

\section{Convergence of Frequency Map Analysis.}

\begin{definition}
A weight function is a positive,   even $C^\infty$ function  $\chi$ on
$[-1,1]$ such that $\frac{1}{2} \int_{-1}^{1} \chi(t) dt = 1$. We will call the
transform of $\chi$, 
 the
$C^\infty$ function on $\bbbr$,  $\vphi_\chi$, defined as   \be
\vphi_\chi(x) = \langle e^{ixt}, 1\rangle_1^\chi  \qquad
\hbox{where}  \qquad
\langle f(t), g(t)\rangle_T^\chi = \Frac{1}{2T}\int_{-T}^T f(t)\bar g(t)
\chi(t/T) \, dt 
\ee
\end{definition}
\begin{lemme}
\label{lemw}
If $\chi$ is a weight function, we have  $\vphi_\chi(0) =1$;  $\vphi'_\chi(0) =
0$; $\vphi''_\chi(0)< 0$.
For all $n\in \bbbn$, we have $\pabs{\vphi_\chi^{(n)}(x)} < 1$, and there
exists $ M_n >0$ such that 
$\pabs{x\vphi^{(n)}(x)}\leq M_n$  on $\bbbr$.  
\end{lemme}
\noindent
{\it Proof.}
As $\chi(t)$ is even,  $\vphi_\chi(x)$ is a real and even function, which implies  $\vphi'_\chi(0) =0$.
We have $\vphi''_\chi(0)< 0$   as it is the integral of a  strictly negative
function.  We have also $\pabs{\vphi^{(n)}(x)} \leq 
\frac{1}{2} \int_{-1}^{1} \abs{t}^n \chi(t) dt <  \frac{1}{2} \int_{-1}^{1}  \chi(t)
dt = 1$. On the other hand, using integration by parts, we obtain
$\pabs{x\vphi^{(n)}(x)} \leq (\chi(1)+\chi(-1))/2 + \frac{1}{2} \int_{-1}^{1}
\chi'(t) dt +n$.
\begin{definition}
Let  $\nu = (\nu_1,\nu_2,\dots,\nu_n)$, and 
\be
\label{eq.ft}
f(t) = e^{i \nu_1 t} + \sum_{k\in \bbbz^n-(1,0,\dots,0)} a_k e^{i \langle k,\nu \rangle t};
\quad a_k\in
\bbbc
\ee
a quasi periodic function on $\bbbr$. In the frequency
map analysis, the approximation $\nu_1^T$ of $\nu_1$  is obtained as the value of $\sigma$ 
for which   
$
\vphi(\sigma) = \abs{ \langle f(t), e^{i\sigma t}\rangle_T^\chi }  
$
is  maximum in a neighborhood of $\nu_1$.
\end{definition}

\begin{theo}
\label{th1}
Let $\chi$ be a weight function, and $\vphi=\vphi_\chi$ its transform with
the asymptotic expressions when $x \rightarrow \infty$
\be
\vphi(x) = \Frac{g_0(x)}{x^n} +o\l( \Frac{1}{x^n}\r),\quad
\vphi'(x) = \Frac{g_1(x)}{x^n} +o\l( \Frac{1}{x^n}\r),\quad
\vphi''(x) = \Frac{g_2(x)}{x^n} +o\l( \Frac{1}{x^n}\r),
\ee
where $n \geq 1$ and where $g_0(x), g_1(x)$ and $g_2(x)$ are bounded on $\bbbr$.
Let $f(t)$ be a quasi periodic function on the form  (\ref{eq.ft}),
and for all $k$, $\Om_k = \langle k,\nu \rangle -\nu_1$; and assume that 
$
\sum_k \pabs{\Frac{a_k}{\Om_k^p}}$  is convergent for $p=0,1$, and $n$.
Then for  $T\rightarrow +\infty$, $\nu_1^T \rightarrow \nu_1$ and 
\be
\nu_1 - \nu_1^T = \Frac{-1}{\vphi''(0) T^{n+1}} \sum_k \Frac{\Re(a_k)}{\Om_k^n}g_1(\Om_k T) 
 + o\l( \Frac{1}{T^{n+1}}\r)
 \llabel{eq.cor}
\ee
\end{theo}

\noindent
{\it Proof.} With   $x= (\nu_1-\s)T$, $\nu_1^T$ will be obtained for    the maximum
value of the modulus of
\be
\psi(x) = \vphi (x)+ \sum_{k\in \bbbz^n-(1,0,\dots,0)}\, a_{k} \vphi \l(x +\Om_k T\r) 
\ee
or in an equivalent way, the maximum of its square
  $\psi(x)\bar\psi(x)$. This maximum thus fulfills the condition
$\psi'(x)\bar\psi(x) + \psi(x)\bar\psi'(x) = 0$,  that is

\be
\EQM{
F(x,T)=\vphi(x)\vphi'(x) &+  \sum_{k,l} \Re(a_k \bar a_l)\ \vphi(x+\Om_k
T)\vphi'(x+\Om_l T) \cr &+ \,  \sum_{k}\Re(a_{k})\l(\vphi(x)\vphi'(x+\Om_k T) +
\vphi'(x)\vphi(x+\Om_k T)\r) = 0
\ .}
\llabel{eq.FxT}
\ee

\begin{lemme}
\label{lem3}
Let $\chi$ be a weight function, and $\vphi = \vphi_\chi$ its 
transform. 
Let us also assume that the series $\sum_k \abs{a_k}$ and 
$\sum_k \abs{\frac{a_k}{\Om_k}}$ are convergent with sum $S_0$ and $S_1$. 
Then for all 
$A>0$, 
\be
\lim_{T\rightarrow +\infty} \sum_{k} a_k \vphi(x+\Om_k T) =0
\ee 
uniformly with respect to $x \in [-A,A]$.
\end{lemme}
\noindent
{\it Proof.} Let ${\cal F}(x,T) = \sum_{k} a_k \vphi(x+\Om_k T)$ and $\eps>0$. With a
Taylor expansion of
$\vphi$ at order $n$, we have
\be
{\cal F}(x,T) =  \sum_{p=0}^n \Frac{x^p}{p!} \sum_{k} a_k \vphi^{(p)}(\Om_k T)
+ \Frac{x^{n+1}}{(n+1)!}\sum_{k} a_k\vphi^{(n+1)}(\xi_k) \quad \hbox{where } \xi_k\in
\bbbr \ .
\ee
We have from lemma \ref{lemw}
\be
\l\vert\Frac{x^{n+1}}{(n+1)!}\sum_{k} a_k\vphi^{(n+1)}(\xi_k)\r\vert \le
\Frac{A^{n+1}}{(n+1)!}S_0
\ee
and this expression can be made arbitrarily small for sufficiently large values of $n$.
On the other hand, with the notations of lemma \ref{lemw}
\be
\l\vert \sum_{k} a_k \vphi^{(p)}(\Om_k T) \r\vert =
 \Frac{1}{T} \abs{\sum_{k} \Frac{a_k}{\Om_k } (\Om_k T)\vphi^{(p)}(\Om_k T)}
\le \Frac{M_p}{T}S_1 \ .
\ee
For any $\eps>0$,  there exists thus $n_0\in  \bbbn$ such that for all 
$n>n_0$, and all $x\in [-A,A]$,
\be
\abs{{\cal F}(x,T)} =  \Frac{S_1}{T}\sum_{p=0}^n \Frac{A^p}{p!}M_p \ + \ \eps
\ee
and for $T$ sufficiently large, the sum will be arbitrarily small, which ends 
the demonstration.

\begin{cor}
With the hypothesis of lemma \ref{lem3}, we have
\be
\lim_{T\rightarrow +\infty} F(x,T) = \vphi(x)\vphi'(x)
\ee
uniformly for $x\in [-A,+A]$. The function $F(x,T)$ can thus be continuously extended 
in a function $F^*(x,T) : [-A,A]\times ]0,+\infty] \rightarrow \bbbr$.
\end{cor}

We have also for the derivative
\be
\EQM{
\Dron{F(x,T)}{x}  &=  \vphi(x)\vphi''(x) +\vphi'^2(x)  \cr
&+ \sum_k \Re(a_k)\l\{ 2\vphi'(x)\vphi'(x+\Om_k T) + 
\vphi(x)\vphi''(x+\Om_k T) + \vphi''(x) \vphi(x+\Om_k T) \r\}  
\cr 
&+   \sum_{k,l} \Re(a_k \bar a_l) \l\{\vphi(x+\Om_k T)\vphi''(x+\Om_l T) +\vphi'(x+\Om_k T)
\vphi'(x+\Om_l T)  \r\} \ .
}
\ee

With the hypothesis of lemma \ref{lem3}, we will have 
\be
\lim_{T\rightarrow  +\infty} \Dron{F(x,T)}{x} =  \vphi(x)\vphi''(x) +\vphi'^2(x) 
\label{dron}
\ee
uniformly for $x\in [-A,+A]$. The extension $F^*(x,T) $ thus has a continuous first derivative
on $[-A,A]\times ]0,+\infty]$. From Lemma \ref{lemw} and (\ref{dron}), we have 
$\partial F^*/\partial x \,(0,+\infty) = \vphi''(0) <0$, we can  apply the
implicit function theorem (Schwartz, 1992, p.176) to $F^*(x,T) $ at
the point 
$(0,+\infty)$. Thus, there exists a neighborhood
$U$ of $+\infty$,   and a unique continuous map 
$\tx(T) = x$  such that
\be
\lim_{T \rightarrow +\infty} \tx(T) =  0 ; 
\qquad F(\tx(T), T) = 0 \quad \hbox{for all\ } T\in U \ .
\llabel{eq.tx}
\ee

Once the existence of  $\tx(T)$ verifying  (\ref{eq.tx}) is known, we can obtain 
a generalized equivalent for $\tx(T)$ when $T \rightarrow +\infty$.

\begin{lemme}
\label{l4}
Let $\vphi(x)$ be a   function on $\bbbr$, and $n\in \bbbn, n>0$, and assume
that for all $0\leq p \leq n$ there exists $M_p >0$ such that 
$\abs{x^p\vphi(x)} \le M_p$ for all $x \in \bbbr$. Then, for all $A>0$,
and all $0\leq p \leq n$, there exists $N_p >0$ such that 
\be
\abs{X^p\vphi(x+X)} \leq N_p \qquad \forall (x,X) \in [-A,+A]\times\bbbr \ .
\ee
\end{lemme}
\noindent
{\it Proof.} We have $X^p\vphi(x+X) = (X+x)^p\vphi(x+X)
-\sum_{k=0}^{p-1} C_p^k X^kx^{p-k}\vphi(x+X)$ and the proof is 
obtained by recurrence  on $p$.

\begin{lemme}
\label{l5}
Let $\vphi(x)$ be a $C^\infty$ function on $\bbbr$, and assume that 
\be
\vphi(x) = \Frac{g_0(x)}{x^n} +o\l( \Frac{1}{x^n}\r),\qquad
\vphi'(x) = \Frac{g_1(x)}{x^n} +o\l( \Frac{1}{x^n}\r),
\ee
where $g_0(x)$ and $g_1(x)$ are bounded on $\bbbr$, and  that the series 
$\sum_k \abs{\frac{a_k}{\Om_k^n}}$ is convergent. Then for $T \rightarrow +\infty$,
\be
\sum_{k} a_k \vphi(x(T)+\Om_k T) = \sum_{k} \Frac{a_k}{\Om_k^nT^n} g_0(\Om_k T) + o\l(
\Frac{1}{T^n}\r)  \ .
\ee 
\end{lemme}
\noindent
{\it Proof.} 
Let
\be
{\cal F}(T) = T^n\l(\sum_{k} a_k \vphi(x(T)+\Om_k T) - 
\sum_{k} \Frac{a_k}{\Om_k^nT^n} g_0(\Om_k T))\r) \ .
\ee
Thus
\be
{\cal F}(T) =  \sum_{k} \Frac{a_k}{\Om_k^n} {\cal F}_k(T)
\label{eq.FT}
\ee
with ${\cal F}_k(T)=(\Om_kT)^n
 \vphi(x(T)+\Om_k T) -g_0(\Om_k T)$. 
From lemma \ref{l4},   $\abs{{\cal F}_k(T)} $ is bounded, and the series
(\ref{eq.FT}) is uniformly convergent with respect to $T$. We thus have
\be
\lim_{T \rightarrow +\infty} {\cal F}(T) = 
\sum_{k} \Frac{a_k}{\Om_k^n} \lim_{T \rightarrow +\infty} {\cal F}_k(T)
\ee
But $\vphi(x(T)+\Om_k T) -\vphi(\Om_k T) = x(T)\vphi'(\xi_k(T))$
with $\abs{\xi_k(T) -\Om_k T}\leq \abs{x(T)}$. Thus, as $x^n\vphi'(x)$ is
bounded on $\bbbr$, we have for all $k$
\be
\EQM{
\lim_{T \rightarrow +\infty} {\cal F}_k(T)= \crm
 \lim_{T \rightarrow +\infty} x(T)\Frac{(\Om_kT)^n}{\xi_k^n(T)}
\xi_k^n(T) \vphi'(\xi_k(T)) +
  (\Om_kT)^n
 \vphi(\Om_k T) -g_0(\Om_k T) = 0 \ ,
 }
\ee
which ends the proof of the lemma.
We can now complete the proof of theorem \ref{th1}.
Assuming the hypothesis of theorem \ref{th1}, we can 
 apply this lemma to $\vphi$ and $\vphi'$. 
We have thus
\be
\EQM{
\sum_{k} a_k \vphi(x(T)+\Om_k T) = \sum_{k} \Frac{a_k}{\Om_k^nT^n} g_0(\Om_k T) + o\l(
\Frac{1}{T^n}\r) \crm
\sum_{k} a_k \vphi'(x(T)+\Om_k T) = \sum_{k} \Frac{a_k}{\Om_k^nT^n} g_1(\Om_k T) + o\l(
\Frac{1}{T^n}\r) \ .
\llabel{eq.33}
}
\ee
We have also $\vphi(x) = 1 +o(x)$,
and $\vphi'(x) = \vphi''(0)x + o(x)$, and  $\lim_{T\rightarrow +\infty} x(T) = 0 $. We can thus make the
expansion of the expressions involved in Eq.\ref{eq.FxT} and obtain, by an expansion at first order in $x$
and order $n$ in $1/T$
\be
\vphi''(0)x(T) +o(x(T)) + \sum_{k} \Frac{\Re(a_k)}{\Om_k^nT^n} g_1(\Om_k T) + o\l(
\Frac{1}{T^n}\r) =0 \ ,
\ee
from which  the proof of the theorem ends easily. Indeed,
\be
T^n x(T) = \Frac{-1}{\vphi''(0) +\Frac{o(x(T))}{x(T)}}\left(\sum_{k} \Frac{\Re(a_k)}{\Om_k^n} g_1(\Om_k T) +
o(1)\right),
\ee
As   $\lim_{T\rightarrow  +\infty} x(T) = 0 $, and $\vphi''(0) \neq 0$, then $T^n x(T)$ is bounded and thus 
$o(x(T))=o(1/T^n)$ from which
\be
 x(T) = -\sum_{k} \Frac{\Re(a_k)}{\vphi''(0)\Om_k^nT^n} g_1(\Om_k T) +o\l(
\Frac{1}{T^n}\r) \ .
\llabel{eq.36}
\ee
This end the proof of theorem 1.
We can then prove the following result 
\begin{theo}  
With the hypothesis of theorem 1,  for $T \rightarrow +\infty$
\be
\nu_1 -\nu_1^T =  \Frac{-1}{T\vphi''(0)} \sum_k \Re(a_k) \vphi'(\Omega_k T) + O(1/T^{2n+1})
\ee
\end{theo}
This improvement requires the following lemma :
\begin{lemme}   
\label{l6p}
Let $\chi$ be a weight function,   $\vphi = \vphi_\chi$ its 
transform, and assume that 
\be
\vphi(x) = \Frac{g_0(x)}{x^n}  +o\l( \Frac{1}{x^{n}}\r),\qquad
\vphi'(x) = \Frac{g_1(x)}{x^n} +o\l( \Frac{1}{x^{n}}\r),\quad
\vphi''(x) = \Frac{g_2(x)}{x^n} +o\l( \Frac{1}{x^{n}}\r),
\llabel{eq.28}
\ee
where $g_0(x),g_1(x),g_2(x)$   are bounded on $\bbbr$, and  that the series 
$\sum_k \abs{\frac{a_k}{\Om_k^p}}$ is convergent for $p=0,n$. Then for $T \rightarrow +\infty$,
\be
\EQM{
\sum_{k} a_k\vphi (x(T)+\Om_k T) = \sum_{k} {a_k} \vphi (\Omega_k T) + O(1/T^{2n})\crm
\sum_{k} a_k\vphi'(x(T)+\Om_k T) = \sum_{k} {a_k} \vphi'(\Omega_k T) + O(1/T^{2n})
}
\ee
\end{lemme}

\noindent
{\it Proof.} We will make the proof of the first relation (for $\vphi(x)$); the argument for 
$\vphi'(x)$ is the same.
Let
\be
\vphi(x(T)+\Om_k T) = \vphi(\Om_k T) + x(T) \vphi'(\Om_k T) +\Frac{x(T)^2}{2}\vphi''(\zeta_k)
\llabel{eq.30}
\ee
where $\zeta_k\in[\Om_k T,\Om_k T+x(T)]$, and 
\be
{\cal F}(T) = T^{2n}\l(\sum_{k} a_k \left[\vphi(x(T)+\Om_k T) - \vphi(\Omega_k T)\right]\r)
\ee
that is,
\be
{\cal F}(T) = T^{2n} \sum_{k} a_k 
\left[x(T) \vphi'(\Om_k T) +\Frac{x(T)^2}{2}\vphi''(\zeta_k)\right] =\sum_{k}   {\cal F}_k(T)
\llabel{eq.32}
\ee

From  theorem 1,   
$\abs{T^{n}x(T)}$ is  strictly bounded by a constant $N$. Moreover, from the hypothesis of the lemma, $x^n \vphi'(x) $ is also bounded by 
a constant $N'$, and thus $\abs{\Om_k^{n}T^n \vphi'(\Om_kT)} < N'$, and
\be
 \abs{a_k  T^{2n} x(T) \vphi'(\Om_kT) } 
 <   N N' \abs{\Frac{a_k}{ \Om_k^{n}}}      
\ee

Finally, $\abs{\vphi''(x)} < 1$ (lemma 1), and as  
$\abs{T^{2n}x(T)^2} < N^2$  
\be
  \abs{T^{2n}\,a_k \Frac{x(T)^2}{2}\vphi''(\zeta_k)}  
<  \Frac{N^2}{2}   \abs{a_k}  \ .
\ee
In  (\ref{eq.32}) $\abs{{\cal F}_k(T)} $ is thus bounded by 
$u_k = NN' \abs{\Frac{a_k}{ \Om_k^{n}}} +  \Frac{N^2}{2}   \abs{a_k}  $
which is the  general term of a convergent series of sum $S$, and 
\be
\abs{{\cal F}(T)} \leq \sum_{k}   \abs{{\cal F}_k(T)} \leq \sum_{k} u_k = S
\ee
which ends the proof of the lemma. The proof of the theorem can now be completed.

\noindent
{\it Proof of theorem 2 :} 
From theorem 1, we already know that  when $T \rightarrow +\infty$, the solution 
$x(T)$ of EQ. (\ref{eq.FxT}) tends to zero and  from (\ref{eq.36}),
\be
x(T) = O(1/T^n)
\ee
We have thus also the expansions  for $T \rightarrow +\infty$
\be
\EQM{
\vphi(x(T)) = 1 + O(1/T^{2n}) \crm
 \vphi'(x(T)) = \vphi''(0) x(T) + O(1/T^{3n})= O(1/T^{n}) ; 
 }
\ee
while Eq.(\ref{eq.33}) gives 
\be
\EQM{
\sum_{k} a_k \vphi(x(T)+\Om_k T) =  O(1/T^{n}) \crm
\sum_{k} a_k \vphi'(x(T)+\Om_k T) = O(1/T^{n}) \ .
}
\ee
Using these expressions, we can now expand 
all expressions in Eq.\ref{eq.FxT} up to $O(1/T^{2n})$ which gives
\be
\vphi''(0) x(T)  + \sum_k \Re(a_k) \vphi'(x(T)+\Omega_k T) = O(1/T^{2n}) 
\ee
and with lemma \ref{l6p}
\be
 x(T)  = -\Frac{1}{\vphi''(0)} \sum_k \Re(a_k) \vphi'(\Omega_k T) + O(1/T^{2n})
\ee
\subsection{Amplitudes}
Once the estimate of the precision on the frequencies is obtained, the precision on hte amplitudes is 
easily calculated. Indeed, with the same notations, the amplitude of the first term of 
$f(t)$ in eq.\ref{eq.ft} (which exact value is 1), will be
\be
a_1^T = < f(t), \e^{i\nu_1^T t}> = \vphi(x(T)) 
+ \sum_k a_k \vphi(x(T)+\Omega_k T) \ .
\ee
Using the above estimates, one then finds that the error on the amplitude is 
\be
a_1^T-1 =  \sum_k a_k \vphi(\Omega_k T)  + O(1/T^{2n}) = O(1/T^{n}) \ .
\ee
\section{Cosine windows}
\begin{figure}[h]
\includegraphics*[scale=0.6]{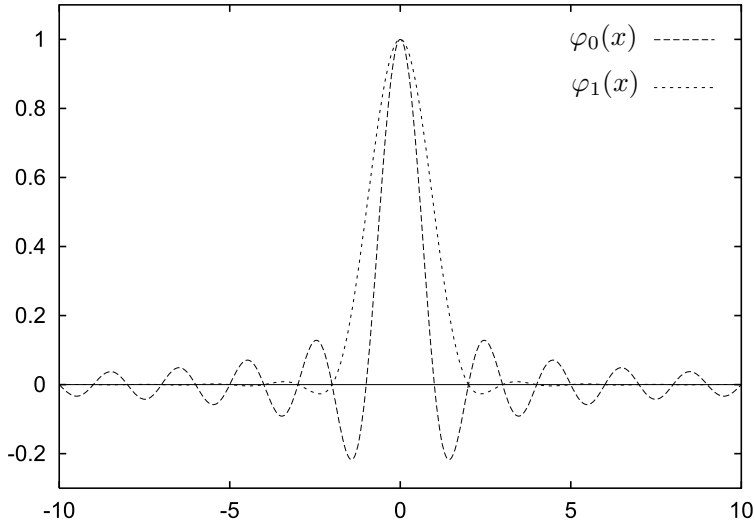}   
\caption{Graph of $\vphi_0(x)$ and $\vphi_1(x)$.} 
\label{porqb.fig1} 
\end{figure}
We can now apply this theorem to specific examples of weight functions. The most simple 
weight function will be defined on $[-1,1]$ as $\chi_0(t) = 1$, with the associate transform 
\be
\vphi_{\chi_0}(x) = \Frac{\sin x}{x}
\ee
With this window (or absence of window), the decreases of the sidelobes of the transform is 
slow, (as $1/x$). This is due to the discontinuity at $-1$ and $+1$ of the 
window function $\tilde \chi_0(x)=1_{[-1,1]}$, considered as a function on $\bbbr$. 
In this case, the transform $\vphi_{\tilde\chi_0}(x)$ is the 
usual Fourier transform of $\tilde \chi_0(x)$,and it is known that the 
decrease at infinity of $\vphi_{\tilde\chi_0}(x)$ depends on the class of regularity of 
$\tilde\chi_0(t)$. More precisely,
\begin{prop}
If $f(t)$ is of class 
$C^p$ on $\bbbr$, absolutely integrable, with all derivatives absolutely integrable, then its Fourier transform 
$\hat f(x)$ satisfies 
\be
\hat f(x) = \int_{-\infty}^{+\infty} f(t) \,\e^{ixt} dt = o(1/x^p) \qquad  
\hbox{when}\qquad  \abs{x} \rightarrow  +\infty
\ee
\end{prop}
We will thus improve the decreasing at infinity of the transform function $\vphi_\chi(x)$ by increasing the class of regularity 
of the extension to $\bbbr$ of the window function $\chi(t)$. This can be done by using the window function 
$\chi_1(t) = 1 + \cos \pi \, t$. Indeed, we have then $\chi_1'(t) = -\pi\sin \pi \, t$,  $\chi_1''(t) = \pi^2\cos \pi \, t$ and thus
$\chi_1(1) = \chi_1(-1)=\chi_1'(1) = \chi_1'(-1)=0$, while $\chi_1''(1) = \chi_1''(-1)=\pi^2$.
The extension $\tilde\chi_1(t)$ is thus $C^1$ and  the decreasing at infinity of the sidelobes 
of the transform $\vphi_{\chi_1}(x)$ will be much stronger (Fig.\ref{porqb.fig1}). More generally,
we can search for a trigonometric polynomial $\chi_p$ such that the extension $\tilde \chi_p(t)$ is of class $C^p$ 
over $\bbbr$. Indeed,

\begin{prop}
For all $p\in \bbbn$, let $\chi_p(t) = \frac{2^p(p!)^2}{(2p)!}(1+\cos \pi t)^p$. Then $\chi_p(t)$ is the 
unique trigonometric polynomial $P(\cos\pi t, \sin\pi t)$ of degree $\leq p$  that is a 
weight function of class $C^{2p-1}$. Its associated transform is
\be
\vphi_p(x) =\Frac{1}{2}\int_{-1}^1 e^{ixt} \ \chi_p(t)\ dt = 
\Frac{(-1)^p\pi^{2p} (p!)^2 \sin x}{x (x^2-\pi^2)\cdots(x^2-p^2\pi^2)} \ .
\ee
\end{prop}

\noindent
{\it Proof.} As a weight function is even, $P$ can be expressed  uniquely as a polynomial in 
$\cos \pi t$. Setting $u(t)=\cos \pi t$, the derivatives of $P(u(t))$ are given by the Faa di Bruno formula
\be
P(u)^{(k)} = \sum_{r=1}^{k} \sum_{k_1+\cdots +k_r=k} \Frac{1}{r!} \Frac{k!}{k_1!\cdots k_r!}
P^{(r)}(u).u^{(k_1)}\cdots u^{(k_r)}
\llabel{eq.54}
\ee
where $k_1, \dots, k_r \geq 1$.
As for $k \geq 0$, 
\be
\EQM{
u^{(2k )}(t)&= \pi^{2k } (-1)^{k}\cos(\pi t) ,\crm 
u^{(2k+1)}(t)&= \pi^{2k+1} (-1)^{k+1}\sin(\pi t) ;
}
\ee
For $k\geq 0$, we have $u^{(2k+1)}(1)=u^{(2k+1)}(-1)=0$, and  $u^{(2k)}(1)=u^{(2k)}(-1)=(-1)^{k+1}$. 
The only index $k_i$ for which the contribution in Eq.\ref{eq.54} is not zero are thus even, and we obtain
\be
\EQM{
P(u)^{(2k+1)}(1)=P(u)^{(2k+1)}(-1) =0 \crm
P(u)^{(2k)}(1)=P(u)^{(2k)}(-1) = \sum_{r=1}^{k} \sum_{k_1+\cdots +k_r=k} \Frac{(-1)^{k+r}}{(r)!} 
\Frac{(2k)!}{(2k_1)!\cdots (2k_r)!}
P^{(r)}(-1) \ .
}
\llabel{eq.56}
\ee
To Assume that the weight function is of class $C^{2p-1}$ is equivalent to assume that 
$P(u)^{(k)}(1)=P(u)^{(k)}(-1) =0$ for $k=0,\dots 2p-1$, that is, from (\ref{eq.56}), 
$P^{(k)}(-1) = 0$, for $k=0,1,\dots,p-1$. P is thus a polynomial of the form 
$P(X) = (1+X)^p\, Q(X)$, where $Q(X) $ is a polynomial, and the  proposition follows.
The constant is determined such that $1/2\int_{-1}^1 \chi_p(t) dt = 1$.
 We have also easily the asymptotic expansions for $x\rightarrow \infty$
\be
\vphi_p^{(n)}(x) = 
\Frac{(-1)^{p}\pi^{2p} (p!)^2 \sin^{(n)} x}{x^{2p+1}} + o\l( \Frac{1}{x^{2p+1}} \r) \ ;
\ee
and
the expansion at the origin 
\be
\vphi_p (x) = 1 +\Frac{\vphi_p''(0)}{2}x^2 + o(x^2)
\quad \hbox{with}\quad 
\vphi_p''(0) = - \Frac{2}{\pi^2}\l(\Frac{\pi^2}{6} -\sum_{k=1}^p\Frac{1}{k^2}\r)
\ee
Thus, if  $f(t)$ is a quasi periodic function  of the form (\ref{eq.ft}),  for which 
$\sum_k \pabs{\Frac{a_k}{\Om_k^m}}$  is convergent for  $m=0,1$, and $2p+1$; 
for  $T\rightarrow +\infty$, 
\be
\nu_1 - \nu_1^T = \Frac{(-1)^{p+1}\pi^{2p} (p!)^2 }{\vphi_p''(0) T^{2p+2}} \sum_k
\Frac{\Re(a_k)}{\Om_k^{2p+1}}\cos(\Om_k T) 
 + o\l( \Frac{1}{T^{2p+2}}\r)
 \llabel{eq.59}
\ee

Using theorem 2, we have as well  for the cosine windows $\vphi_p(x)$, 
\be
\nu_1 -\nu_1^T = \Frac{-1}{T\vphi_p''(0)} \sum_k \Re(a_k) \vphi_p'(\Omega_k T) + o(1/T^{4p+2})
 \llabel{eq.60}
\ee

The computation of the derivative $\vphi_p(x)$ 
is  easy if we write $ \vphi_p(x) = C \sin x / D(x) $ 
where 
\be
C =  (-1)^p\pi^{2p} (p!)^2 \ ; \qquad  D(x) = x (x^2-\pi^2)\cdots(x^2-p^2\pi^2)
\ee
Then 
\be
\vphi_p'(x) = \Frac{C}{D(x)} \left[ \cos x -\sin x \Frac{D'(x)}{D(x)}\right] 
\ee
with 
\be
 \Frac{D'(x)}{D(x)} =\Frac{1}{x} + \Frac{2x}{x^2 -\pi^2 } + \cdots   + \Frac{2x}{x^2 -p^2\pi^2 }
\ee

\subsection{Exponential window}
\label{expw}
The cosine weight function of order $p$, $\chi_p(t)$ provides a window function 
of class $C^{2p-1}$, and thus improves the convergence of 
the frequency map algorithm (theorem 1 and 2). Although it will not be always very useful in 
practice, 
we can  also  consider $C^{\infty}$ windows functions,  for which 
$\chi^{(n)}(-1)=\chi^{(n)}(1)=0$ for all $n$, as for example $\chi^*(t) = c\exp(-1/(1-t^2))$, with 
$c\approx 0.22199690808403971891$. In
this case, we have for all values of $n,m$, $\lim_{x\rightarrow +\infty} x^n\vphi_\chi^{(m)}(x)
=0$. Thus, if $\sum_k \pabs{\frac{a_k}{\Om_k^n}}$
are convergent for all $m$, as for a KAM solution (see below), we will have for all $n$, 
$\lim_{T\rightarrow +\infty} T^n(\nu_1-\nu_1^T) =0$. Let us note that for this exponential window,
$\vphi_\chi''(0)\approx -0.035100738376487704994$.

\subsection{KAM solutions}
If $f(t) $ is a KAM solution, with rotation vector $(\nu)$, that fulfills a   
a diophantine condition (\ref{eq.3}),
it can be expressed as a quasiperiodic series 
of the form (\ref{eq.5}),
\be
f(t) =\sum_k a_k \e^{i < k,\theta >} \qquad \hbox{with} \qquad \theta_i=\nu_i t
\llabel{eq.64}
\ee
 that is analytical with respect to the angles $(\theta)$.
 The series $\sum_k \pabs{\frac{a_k}{\Om_k^p}}$
are then convergent for all $p$. 
Indeed, we have then
\be
\abs{\frac{a_k}{\Om_k^p}} \leq \Frac{\abs{a_k}}{\kappa_\eps^p}(1+\abs{k})^{pm} \ .
\ee
But due to the analyticity of the series $(\ref{eq.64})$ with respect to $\theta$, there exists 
a small $s > 0$ such that 
\be
\abs{a_k} \leq \norm{f} \e^{-\abs{k}s}
\ee
As   $\abs{k} = \abs{k_1}+\cdots +\abs{{k_n}}$, where $n$ is the number of degrees of freedom 
of the Hamiltonian system, 
\be
(1+\abs{k}) \leq (1+\abs{k_1})(1+\abs{k_2})\dots (1+\abs{k_n})
\ee
\be
\EQM{
\sum_k \abs{\frac{a_k}{\Om_k^p}} &\leq 
\sum_{k_1,\dots,k_n} \Frac{\norm{f} }{\kappa_\eps^p}(1+\abs{k_1})^{pm} \dots (1+\abs{k_n})^{pm} \e^{-\abs{k}s} \crm
&= \Frac{\norm{f} }{\kappa_\eps^p}\prod_{i=1,\dots,n} \sum_{k_i\geq 0}(1+\abs{k_i})^{pm}\e^{-\abs{k_i}s} \ .
}
\ee
These latest series are convergent which ends the proof.
In fact, this case will be the most useful, as the solutions we are looking for 
will usually be the quasiperiodic KAM solutions of an analytical Hamiltonian system.
In this case, 
 using the weight $\chi_p$, we
can obtain an accuracy for  the frequency analysis  of order $1/T^{2p+2}$ on these solutions,
which gives in particular $1/T^4$ when using the  Hanning window ($\chi_1$).

\section{Quasiperiodic decomposition}
The previous results tell us that with the proposed algorithm, we will recover in 
a very efficient manner the frequency $\nu_1$ of the leading periodic term of 
the quasiperiodic decomposition of
\be
f(t) = \sum_{k} A_k \e^{i \nu_k \, t}
\ee  

A first approximation of the corresponding amplitude $A_1$ will be given by 
\be
A_1^{(1)}= <f, e_1>  \qquad \hbox{where} \qquad  e_1=\e^{i \nu_1\, t}
\ee
we are then left with the remainder
\be
f_1=f - <f,e_1> e_1
\ee
and we can start the process again in order to search for the frequency $\nu_2$ of the next 
term. A technical complexity appears here, 
as $e_1 = \e^{i \nu_1\, t}$ and  $e_2 = \e^{i \nu_2\, t}$ are not orthogonal for $<,>_T^\chi$.
A classical way to overcome this difficulty is given by Gramm-Schmidt orthogonalization process.
Indeed, if $(e_1,\dots, e_n)$ is an independent family of unit vectors, one can construct 
an orthogonal family $(e'_1,\dots, e'_n)$, such that for all $1\leq k\leq n$,
 $Vec(e_1,\dots , e_k) = Vec(e'_1,\dots , e'_k)$, where $Vec(u_1,\dots , u_k)$ denotes the  
vector space generated by   $(u_1,\dots , u_k)$. Indeed,
\be
e'_n =\tilde e_n/\norm{\tilde e_n} \quad \hbox{where}\quad  \tilde e_n= e_n - \sum_{k=1}^{n-1} <e_n,e'_k>e'k \ ,
\ee
and 
\be
f_n=f_{n-1} - <f_{n-1},e'_n> e'_n.
\ee
As $f_{n-1}\perp Vec(e_1,\dots , e_{n-1}) = Vec(e'_1,\dots , e'_{n-1})$, it is interesting to note, for practical computation  that 
\be
<f_{n-1},e'_n> = <f_{n-1},e_n> \ .
\ee
\section{Numerical simulations}

The asymptotic convergence of $(\nu_1-\nu_1^T)$ when $T\rightarrow +\infty$ provides 
a good indication of the possibilities of the method, but in practice, this 
asymptotic behavior will be also limited  by the value of the 
involved constants, and  by numerical accuracy. 
In the following, the previous asymptotic expressions are tested for a very simple example of 
a quasiperiodic function. Instead of using the output of a numerical integration, we prefer  here 
to take directly a quasi periodic function with 2   independent frequencies.
The   tested function is 

\be
F_1(t) = \Frac{1}{1 +1/2 \e^{it} + 1/4\e^{-i \omega t}}
\llabel{eq.77}
\ee
with $\omega = 2.02$. This function is chosen because of the slow 
convergence of its Fourier expansion. The  results are displayed  in figures 2--6 and Table 1
  with a cosine window $\chi_p$ of various order $p=0,1,2,3,4,5$, and an exponential window. 
  A second test function $F_2(t)$ is also constructed with the first 50 periodic terms of 
  the frequency decomposition of $F_1(t)$. 

\subsection{Corrections}
Theorems 1 and 2 provide equivalents of the frequency 
error $(\nu_1-\nu_1^T)$ when $T\rightarrow +\infty$. If a quasiperiodic approximation 
of $f(t)$ is known, for example by a first frequency  analysis, then these estimates can be used 
as  corrections for  frequency. 
It should be noted that when the correction is searched, one 
does not know the exact values of the perturbing frequencies $\Omega_k $ involved in 
the formulas of theorem 1 and 2, but an approximate value $\Omega_k'$. 
Under the hypothesis of theorem 1, one has
\be
\Omega_k' = \Omega_k + O\left(\Frac{1}{T^{n+1}}\right)
\ee
and the expression (\ref{eq.cor}) is as well valid with $\Omega'_k$ instead of $\Omega_k$.
In the same way, 
\be
\vphi''(x) = O\left(\Frac{1}{x^{n}}\right) \quad \hbox{for} \quad x \longrightarrow +\infty,
\ee
\be
\vphi'(\Omega'_k T) = \vphi'(\Omega_k T) + O\left(\Frac{1}{T^{2n}}\right)
\ee
and the estimation of theorem 2 is also valid with $\Omega_k'$ instead of $\Omega_k$. This allows 
to simplify greatly the use of these corrections as the first estimation is sufficient 
in general.

\subsection{Discussion}
\begin{figure}[]
\includegraphics*[scale=0.6]{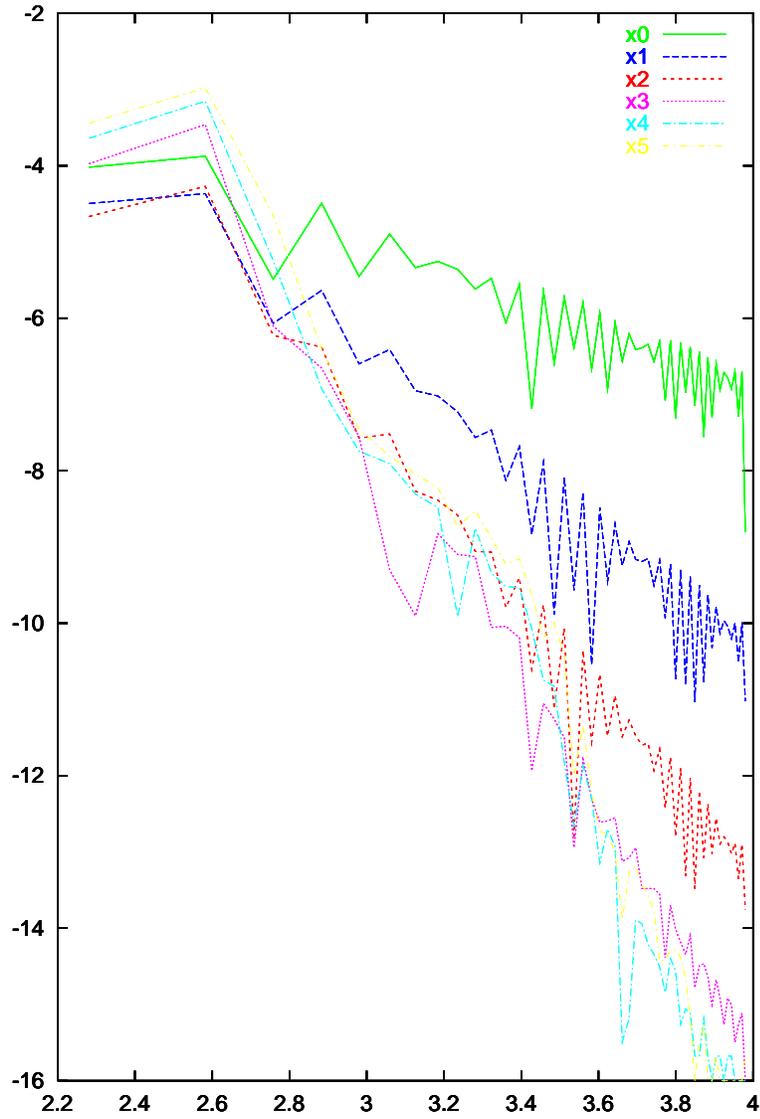}   
\caption{Error on the measured frequency versus duration of the integration
for the first frequency of $F_1(t)$. $\log_{10}(err)$ is plotted versus 
$\log_{10}(T/2 \pi)$. $x_p$ is the result for a window of order $p$ . }
\label{porq.fig2} 
\end{figure}

We have first analyzed the  results of the frequency analysis of $F_1(t)$ for different 
windows $\chi_p$, and the results of the precision of the determination 
of the frequency versus the length of the time interval $T$ are given in 
Figure 2. 
For each case, the results are fitted with a line 
and the the results are gathered in  Table 1
where $a_0$ (column 5) is the slope  of the straight line least squares fitted   to the curves of 
Fig.\ref{porq.fig2}.
The same study is then repeated with a correction step provided by theorem 1 and formula 
(\ref{eq.59}) using 10 and 50 terms . The results are given  in Figures 3--4 and 
in column $a_{10}$ and $a_{50}$ of Table 1. These values need to be compared with 
the theoretical value obtained with the full correction of theorem 1, 
given in column $a_\infty$.
The correction obtained with  theorem 2 (formula (\ref{eq.60})) 
is then given in column $a'_{50}$, with a correction 
computed  also with 50 terms.

It should be noted that the correction step  improves the determination 
of the frequency, but only to a certain limit. 
With $F_1(t)$, it is usually hopeless to expect the approximation in $O(1/T^{2n+1})$ given in theorem 2
(column $a'_\infty$ of Table 1).
Indeed, because of  the slow decreases of the amplitude of the quasiperiodic expansion of 
$F_1(t)$, the contribution of the neglected terms will dominate 
beyond a certain level of precision.

We have then performed the same experience with $F_2(t)$. 
In this case, the corrections performed with 50 terms can be considered as optimal, as 
the correction step is now performed with a very precise estimate 
of the remaining  terms.
Indeed, the results that are plotted in Figs  6,7, and gathered in 
columns $b_0, b_{50}, b'_{50}$ of Table 1 are now in much better agreement 
respectively with the  theoretical coefficients $a,a_\infty,a'_\infty$.

It is clear from these numerical simulations that  even for a perfectly 
quasiperiodic function, there are some limitation on  practical  computation of the frequencies,
especially when the quasiperiodic expansion is  slowly converging. 
It is nevertheless striking to 
see that even in the case of a  slow convergence, the agreement  of $a$ and $a_0$ is excellent.
When the correction is used, there is still a good agreement although not as good,
as the  effect of the neglected terms will dominate after a certain level
of precision, unless, as with $F_2(t)$, the quasiperiodic expansion of the 
considered function is fully recovered with good accuracy.

\begin{table}[h]
\begin{tabular}{|c || c|c|c||c|c|c| c||c| c|c| }
\trait
$p $& $-a$ & $-a_{\infty}$ & $-a'_{\infty}$ & $-a_{0}$ &$-a_{10}$ & $-a_{50}$  & $-a'_{50}$& $-b_0$& $-b_{50}$& $-b'_{50}$\\
\trait
0 & 2  &3 & 3 &1.98 & 2.11 & 2.97 &  2.94 & 1.91& 3.15  & 3.11\\
1 & 4  &5 & 7 &3.94 & 4.25 & 4.90 &  4.43 & 3.98& 4.81  & 6.04\\
2 & 6  &7 & 11& 5.82 & 6.45 & 6.35 &  6.11 & 5.89& 7.01 &12.68\\
3 & 8  &9 & 15 &7.64 & 7.99 & 8.01 &  7.87 & 8.69& 9.49 &18.85\\
4 & 10 &11& 19 &8.89 & 9.01 & 9.03 &  8.71 & 11.42&12.48&23.49\\
5 & 12 &13& 23 &9.15 & 9.33 & 9.33 &  8.87 & 13.50&14.47&  26.24\\
$*$ &    &  &    &5.46 &      &      &       & 5.43  &    &      \\
\trait
\end{tabular}
\caption{Slope of the fitted errors on the frequency. $p$ is the order of the window, $a$ the 
exponent of the error from eq. (\ref{eq.59}), $a_0$ the fitted slope with no correction, 
$a_{10}$  and  $a_{50}$  the correction using 10 or 50 periodic terms 
and $a_{\infty} $ the theoretical slope when all terms 
of equation (\ref{eq.59})are used for the correction 
(providing the series is converging) for the determination of the first frequency of $F_1(t)$ (eq.\ref{eq.77}). 
$a'_{\infty} $ is the same with equation (\ref{eq.60}).
$b_0, b_{50},b'_{50}$ are the corresponding value for the function 
$F_2(t)$ obtained as a truncated expansion of $F_1(t)$ with 50 periodic terms. }
\end{table}

Overall, it appears that when the quasiperiodic expansion 
is not fully recovered (up to a certain precision), 
the correction obtained  through equations (\ref{eq.59}) and 
(\ref{eq.60}) does not improve very much the results that can be deduced from
the use of a window of higher order (Table 1). 
Considering that the computation of many terms is   time consuming, and that it is 
not so easy to obtain always the same number of relevant terms for various initial conditions 
corresponding to  different regularity of the solution, I would rather recommend to 
use different orders for the window of the data. The more the solutions are regular, 
the higher one can choose the order of the window. 
Practically, this order is still limited to $p=3,5$ for quasi periodic function, 
and maybe  lower for numerical  solutions of dynamical systems.
When the motion is  more chaotic, $p=1$, or even sometimes $p=0$ may be preferred.
In practice, one would increase the order of the window until the precision seems 
to decreases and use the highest possible value for $p$.

\subsection{Exponential window}
We have also tested the exponential window of section (\ref{expw}), and the results are reported 
for $p=*$ in Table 1, but as it can be seen in Table 1, this was not very successful, compared to 
the cosine windows. Investigation of different exponential windows  should 
probably be continued, as it is not clear why such windows could not 
attain the accuracy reached by the cosine windows $\chi_p$.

\begin{figure}[]
\includegraphics*[scale=0.6]{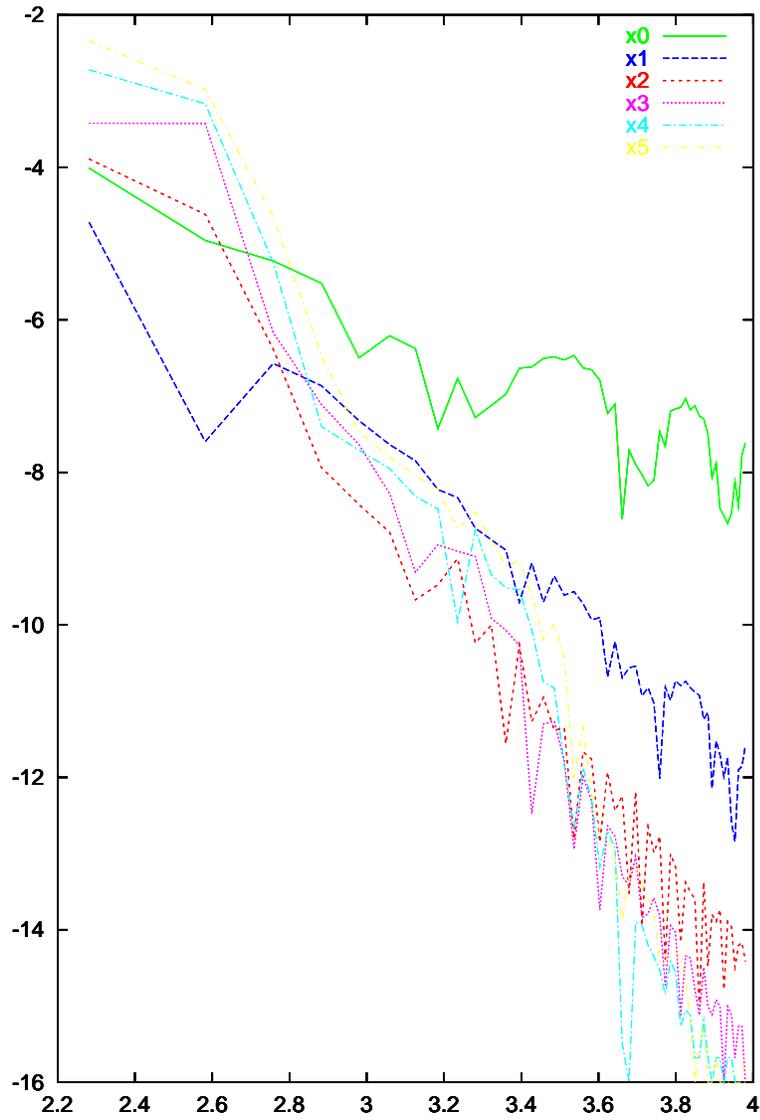}  
\caption{Error on the measured frequency versus duration of the integration
for the first frequency of $F_1(t)$. $\log_{10}(err)$ is plotted versus 
$\log_{10}(T/2 \pi)$. $x_p$ is the result for a window of order $p$ with a correction using 
equation (\ref{eq.59}) with 10 terms.}
\label{porq.fig3} 
\end{figure}

\begin{figure}[]
\includegraphics*[scale=0.6]{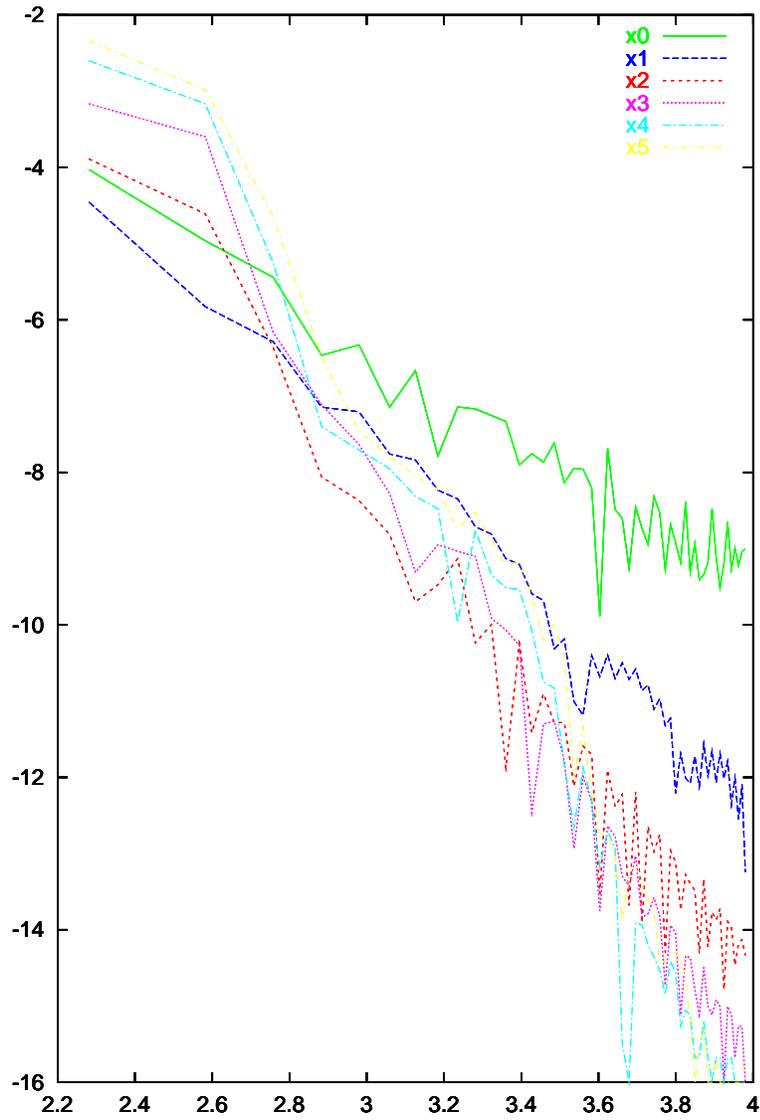}   
\caption{Error on the measured frequency versus duration of the integration
for the first frequency of $F_1(t)$. $\log_{10}(err)$ is plotted versus 
$\log_{10}(T/2 \pi)$. $x_p$ is the result for a window of order $p$ with a correction using 
equation (\ref{eq.59}) with 50 terms.}
\label{porq.fig4} 
\end{figure}

\begin{figure}[]
\includegraphics*[scale=0.6]{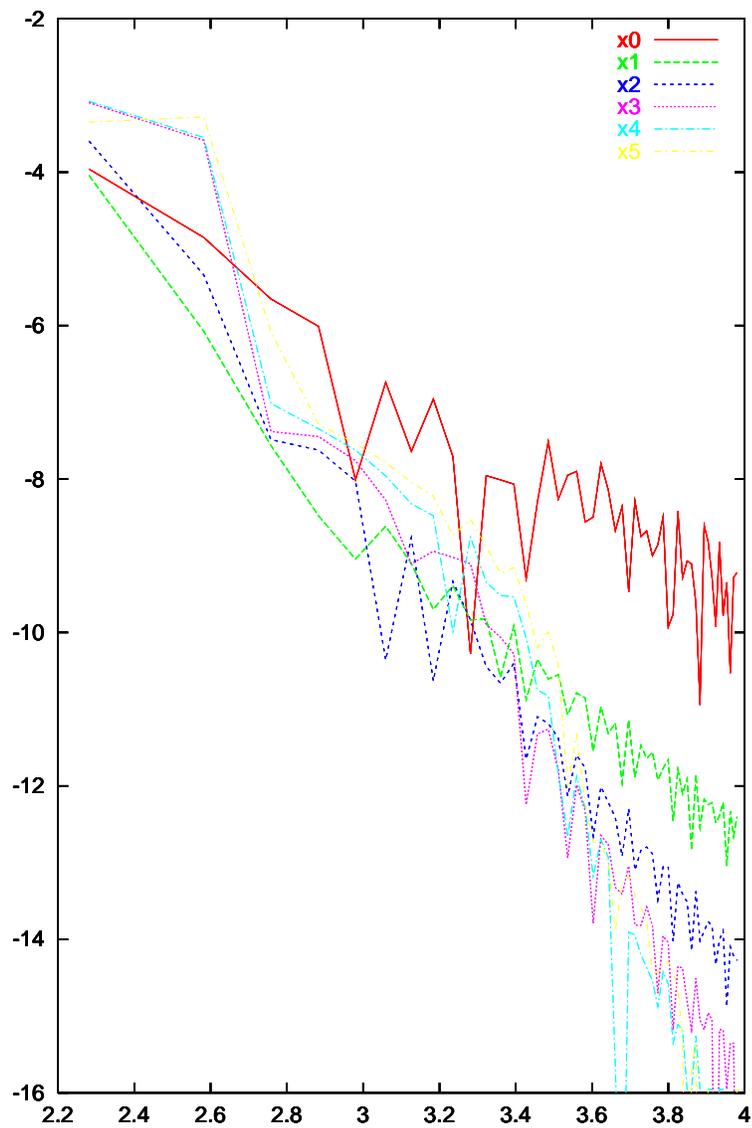}   
\caption{ Error on the measured frequency versus duration of the integration
for the first frequency of $F_1(t)$. $\log_{10}(err)$ is plotted versus 
$\log_{10}(T/2 \pi)$. $x_p$ is the result for a window of order $p$ with a correction using 
equation (\ref{eq.60}) with 50 terms.}
\label{porq.fig5} 
\end{figure}

\begin{figure}[]
\includegraphics*[scale=0.6]{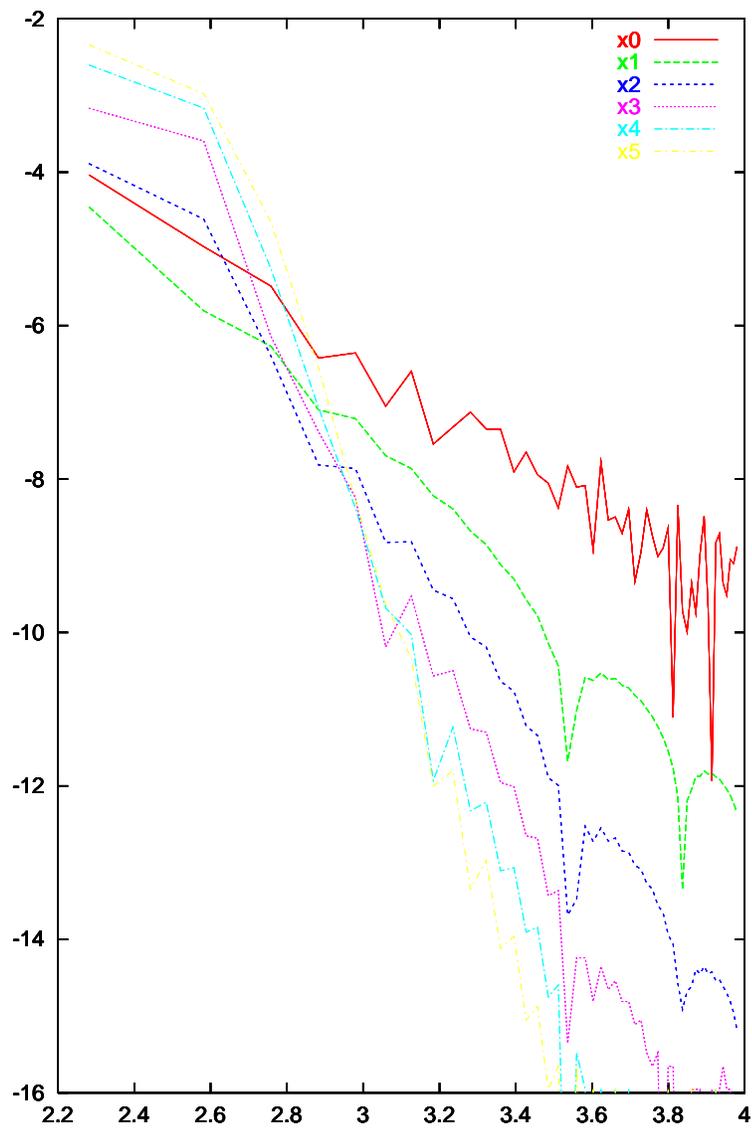}   
\caption{ $F_2(t)$. Error on the measured frequency versus duration of the integration
for the first frequency of $F_2(t)$. $\log_{10}(err)$ is plotted versus 
$\log_{10}(T/2 \pi)$. $x_p$ is the result for a window of order $p$ with a correction using 
equation (\ref{eq.59}) with 50 terms.}
\label{porq.fig5} 
\end{figure}

\begin{figure}[]
\includegraphics*[scale=0.6]{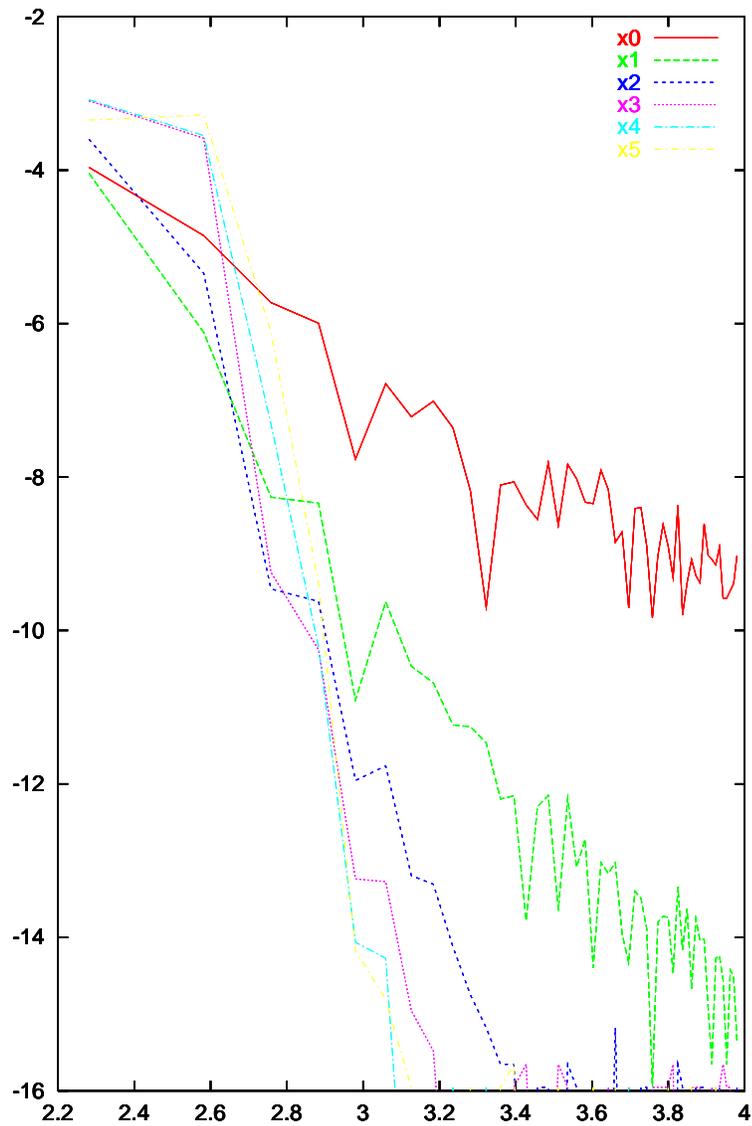}   
\caption{ $F_2(t)$. Error on the measured frequency versus duration of the integration. $\log_{10}(err)$ is plotted versus 
$\log_{10}(T/2 \pi)$ for the first frequency of $F_2(t)$. $x_p$ is the result for a window of order $p$ with a correction using 
equation (\ref{eq.60}) with 50 terms.}
\label{porq.fig5} 
\end{figure}

\clearpage
\section{Beyond Nyquist frequency}

\subsection{Nyquist frequency}

A typical limitation in spectral analysis is given by the so called Nyquist
frequency. Roughly speaking, it will not be possible to determine 
a frequency $\nu$ larger than $\pi/h$, where $h$ is the sampling
time interval of the observed data.
If   $\abs{\nu_0} > \pi /h$, the observed frequency will become $\nu = \nu_0 + k2\pi/h$,
where $k\in \bbbz$ such that $\abs{\nu} < \pi/h$.
This is actually what can be observed using the NAFF algorithm 
detailed above. 
\begin{table}[h]
\begin{tabular}{| r|  r| c|}
\trait
$\nu_0/\pi$ & $\nu/\pi$ & $(\nu_0-\nu)/\pi$  \\ 
\trait
 $     0.990$ & $     0.990$ & $0.000000000000000000$ \\
 $     0.991$ & $     0.991$ & $0.000000000000000000$ \\
 $     0.992$ & $     0.992$ & $0.000000000000000000$ \\
 $     0.993$ & $     0.993$ & $0.000000000000000000$ \\
 $     0.994$ & $     0.994$ & $0.000000000000000000$ \\
 $     0.995$ & $     0.995$ & $0.000000000000000000$ \\
 $     0.996$ & $     0.996$ & $0.000000000000000000$ \\
 $     0.997$ & $     0.997$ & $0.000000000000000000$ \\
 $     0.998$ & $     0.998$ & $0.000000000000000000$ \\
 $     0.999$ & $     0.999$ & $0.000000000000000000$ \\
 $     1.000$ & $    -1.000$ & $2.000000000000000000$ \\
 $     1.001$ & $    -0.999$ & $2.000000000000000000$ \\
 $     1.002$ & $    -0.998$ & $2.000000000000000000$ \\
 $     1.003$ & $    -0.997$ & $2.000000000000000000$ \\
 $     1.004$ & $    -0.996$ & $2.000000000000000000$ \\
 $     1.005$ & $    -0.995$ & $2.000000000000000000$ \\
 $     1.006$ & $    -0.994$ & $2.000000000000000000$ \\
 $     1.007$ & $    -0.993$ & $2.000000000000000000$ \\
 $     1.008$ & $    -0.992$ & $2.000000000000000000$ \\
 $     1.009$ & $    -0.991$ & $2.000000000000000000$ \\
\trait
\end{tabular}
\caption{The function $\exp(i\nu_0 t)$ is evaluated for $t=-1000,1000$, with 
a stepsize $h=1$. For $\nu_0/\pi < h$, the recovered frequency ($\nu$) is recovered 
by NAFF up to the machine precision, but for $\nu_0/\pi > h$, the recovered frequency 
is here $\nu_0-2\pi/h$. }
\llabel{tab2}
\end{table}

\subsection{Multiple timescales problem}
The Nyquist aliasing constraint 
 means that to recover a given period, one needs  to sample the data with at least two points 
per period. On the opposite, in order to determine precisely the long periods,
one needs  that the total interval length $T$ is several time larger than these   periods,
in order to reach the asymptotic rates of theorems 1 and 2, or to be able to 
separate properly close frequencies.

One thus realizes that if good low frequency determination imposes that $T$ is large, 
while high frequencies impose that $h$ is small, we will face a problem  when two 
very different time scales are present in the system.
This is actually the case for planetary systems where short periods (of the order of the year) 
are present, as well as long secular periods, ranging, in the Solar system for example, 
 from 40 000 years  to a few millions of years (Laskar, 1990).

The first frequency analysis of the solar system solutions were made uniquely on the 
output of the secular equations, where only the long periods were present
(Laskar, 1988), but if we want now to perform a frequency analysis 
of the direct output of a numerical integration of Newton's equations, 
without filtering or averaging, both time scales will be present.

For the Solar  system, for example, the determination of the long 
secular frequencies with a good accuracy will require that $T$ is larger than 20 
millions of years, while the sampling $h$ must be smaller than half of the shortest period 
of the system, that is a few days if we consider the Moon motion, while 
$h=5000$ years was enough for the secular equations.
The amount of data to handle becomes then considerable as well as the 
related numerical computations.

\subsection{Multiple timescales solution}
The solution  that can be use to overcome  this problem is simply to 
sample the data with two different sampling intervals that are 
very close $h$ and $h'=h+\eps$.
For a real $x$, we will denote $[x]$ the integer such that 
\be
-\Frac{1}{2} < x-[x] \leq \Frac{1}{2} \ .
\ee
With this notation, the frequency analysis of   $f(t)=\exp (i 2\pi\nu_0\, t)$ will 
give a frequency $2\pi\nu$  such that 
\be
\nu_0 h = \nu  h + k \qquad  \hbox{with} \qquad 
k = [ \nu_0 h ] \ .
\label{eq.84}
\ee
The Nyquist condition is then expressed by the fact that 
$\abs{{\nu_0 h}} < 1/2 $ implies $k=0$.
We assume now that the second sampling time interval $h'=h+\eps$ is such that 
\be
\abs{{\nu_0 \eps}} <  \Frac{1}{2} \ .
\ee
We have then 
\be
\EQM{
\nu_0 h'  = \nu_0(h+\eps) &=\nu h +\eps \nu_0 + k\cr
&= \nu'(h+\eps) + k'
}
\llabel{eq.76}
\ee
with 
\be
k'=k+[\nu h +\eps \nu_0]
\ee
and thus 
\be
\nu'h'-\nu h = \eps \nu_0 - [\nu h +\eps \nu_0] \ .
\label{eq.88}
\ee

As $\abs{ \nu_0\eps} <  1/2$,

\be
[\nu'h'-\nu h] = - [\nu h +\eps \nu_0] \ ,
\ee
and thus, from (\ref{eq.88}),
\be
 \nu_0 =\Frac{1}{h'-h}\left( \nu'h'-\nu h - [\nu'h'-\nu h] \right)\ ,
\ee
which allows to revover in all cases the  true value of the frequencies. Practically,
although this formula can be used, there is an easy way to reduce the  
numerical errors of these computations. Indeed, it is easy to see that 
\be
k = \Frac{h}{h'-h}((\nu'-\nu)h' - [\nu'h'-\nu h]) \ ,
\label{eq.91}
\ee
which allows to recover $\nu_0$ through  (\ref{eq.84}). The main advantage 
of this intermediary step is that as $k$ should be an integer, the computed 
value  can be rounded to its nearest integer
and
\be
\nu_0  = \nu   + \Frac{[k]}{h} \ .
\label{eq.92}
\ee

\begin{table}
\begin{tabular}{| r|  r| r|r|}
\trait
$\nu_0/\pi$ & $\nu/\pi$ & $\nu'/\pi$ & $\nu_f/\pi$ \\ 
\trait
     0.500 &     0.500000 &     0.500500 &  0.500000000000  \\
     1.000 &    -1.000000 &    -0.999000 &  1.000000000000  \\
     1.500 &    -0.500000 &    -0.498500 &  1.500000000000  \\
     2.000 &    -0.000000 &     0.002000 &  2.000000000000  \\
     2.500 &     0.500000 &     0.502500 &  2.500000000000  \\
     3.000 &    -1.000000 &    -0.997000 &  3.000000000000  \\
     3.500 &    -0.500000 &    -0.496500 &  3.500000000000  \\
     4.000 &    -0.000000 &     0.004000 &  4.000000000000  \\
     4.500 &     0.500000 &     0.504500 &  4.500000000000  \\
     5.000 &    -1.000000 &    -0.995000 &  5.000000000000  \\
\trait
   990.000 &     0.000000 &     0.990000 & 990.000000000000  \\
   990.500 &     0.500000 &    -0.509500 & 990.500000000000  \\
   991.000 &    -1.000000 &    -0.009000 & 991.000000000000  \\
   991.500 &    -0.500000 &     0.491500 & 991.500000000000  \\
   992.000 &    -0.000000 &     0.992000 & 992.000000000000  \\
   992.500 &     0.500000 &    -0.507500 & 992.500000000000  \\
   993.000 &    -1.000000 &    -0.007000 & 993.000000000000  \\
   993.500 &    -0.500000 &     0.493500 & 993.500000000000  \\
   994.000 &    -0.000000 &     0.994000 & 994.000000000000  \\
   994.500 &     0.500000 &    -0.505500 & 994.500000000000  \\
   995.000 &    -1.000000 &    -0.005000 & 995.000000000000  \\
   995.500 &    -0.500000 &     0.495500 & 995.500000000000  \\ 
   996.000 &    -0.000000 &     0.996000 & 996.000000000000  \\
   996.500 &     0.500000 &    -0.503500 & 996.500000000000 \\
   997.000 &    -1.000000 &    -0.003000 & 997.000000000000 \\
   997.500 &    -0.500000 &     0.497500 & 997.500000000000 \\
   998.000 &    -0.000000 &     0.998000 & 998.000000000000 \\
   998.500 &     0.500000 &    -0.501500 & 998.500000000000 \\
   999.000 &    -1.000000 &    -0.001000 & 999.000000000000 \\
   999.500 &    -0.500000 &     0.499500 & 999.500000000000 \\
  1000.000 &    -0.000000 &    -1.000000 & 1000.000000000000 \\
\trait
  1000.500 &     0.500000 &    -0.499500 & -999.500000000000 \\
  1001.000 &    -1.000000 &     0.001000 & -999.000000000000 \\
  1001.500 &    -0.500000 &     0.501500 & -998.500000000000  \\
  1002.000 &    -0.000000 &    -0.998000 & -998.000000000000  \\
  1002.500 &     0.500000 &    -0.497500 & -997.500000000000  \\
  1003.000 &    -1.000000 &     0.003000 & -997.000000000000  \\
\trait
\end{tabular}
\caption{The function $\exp(i\nu_0 t)$ is evaluated for $t=-1000,1000$, with 
a stepsize $h=1$, and with  $ h'=1.001$. The reconstruction formula
(\ref{eq.91},\ref{eq.92}) allows to obtain the correct frequency up to $\nu_0/\pi = 1000$.}
\llabel{tab3}
\end{table}

\begin{table}
\small
\thispagestyle{empty}
\begin{tabular}{|  r| r|r|r| r|  r| r|r|r|}
\trait
                        $\nu_{0i}$ &            $P_i$ (years) &$k$& $\abs{A_i}\times 10^{8}$ &   $k_{1i}$&  $k_{2i}$&  $k_{3i}$&  $k_{4i}$&  $k_{5i}$ \\ 
\trait
  $            4.027808$& $  321763.1$& $   0$& $        4412046$ &     $  0$& $  0$& $  1$& $  0$& $  0$\\
  $           28.013748$& $   46263.0$& $   0$& $        1592683$ &     $  0$& $  0$& $  0$& $  1$& $  0$\\
  $       -21264.867867$& $     -60.9$& $  -3$& $          64449$ &     $ -1$& $  2$& $  0$& $  0$& $  0$\\
  $           51.999689$& $   24923.2$& $   0$& $          62876$ &     $  0$& $  0$& $ -1$& $  2$& $  0$\\
  $         1410.142247$& $     919.1$& $   0$& $          38666$ &     $ -2$& $  5$& $  0$& $ -2$& $  0$\\
  $        22703.023863$& $      57.1$& $   4$& $          13141$ &     $ -1$& $  3$& $  0$& $ -1$& $  0$\\
  $       -86525.641213$& $     -15.0$& $ -13$& $          10467$ &     $ -2$& $  3$& $  0$& $  0$& $  0$\\
  $         1386.156306$& $     935.0$& $   0$& $           9940$ &     $ -2$& $  5$& $  1$& $ -3$& $  0$\\
  $        43995.905478$& $      29.5$& $   7$& $           8056$ &     $  0$& $  1$& $  0$& $  0$& $  0$\\
  $       -42557.749483$& $     -30.5$& $  -7$& $           6449$ &     $ -2$& $  4$& $  0$& $ -1$& $  0$\\
  $       -21240.881927$& $     -61.0$& $  -3$& $           4605$ &     $ -1$& $  2$& $ -1$& $  1$& $  0$\\
  $       -21288.853808$& $     -60.9$& $  -3$& $           4294$ &     $ -1$& $  2$& $  1$& $ -1$& $  0$\\
  $      -151786.414559$& $      -8.5$& $ -23$& $           3662$ &     $ -3$& $  4$& $  0$& $  0$& $  0$\\
  $           75.985629$& $   17055.9$& $   0$& $           3493$ &     $  0$& $  0$& $ -2$& $  3$& $  0$\\
  $       109256.678825$& $      11.9$& $  17$& $           3448$ &     $  1$& $  0$& $  0$& $  0$& $  0$\\
  $         1434.128187$& $     903.7$& $   0$& $           2485$ &     $ -2$& $  5$& $ -1$& $ -1$& $  0$\\
  $      -107818.522829$& $     -12.0$& $ -17$& $           2015$ &     $ -3$& $  5$& $  0$& $ -1$& $  0$\\
  $          -19.958132$& $  -64935.9$& $   0$& $           1929$ &     $  0$& $  0$& $  2$& $ -1$& $  0$\\
  $        22679.037922$& $      57.1$& $   3$& $           1852$ &     $ -1$& $  3$& $  1$& $ -2$& $  0$\\
  $         1362.170366$& $     951.4$& $   0$& $           1824$ &     $ -2$& $  5$& $  2$& $ -4$& $  0$\\
  $      -217047.187912$& $      -6.0$& $ -33$& $           1503$ &     $ -4$& $  5$& $  0$& $  0$& $  0$\\
  $       239778.225515$& $       5.4$& $  37$& $           1473$ &     $  3$& $ -2$& $  0$& $  0$& $  0$\\
  $       -86501.655273$& $     -15.0$& $ -13$& $           1217$ &     $ -2$& $  3$& $ -1$& $  1$& $  0$\\
  $          -80.092473$& $  -16181.3$& $   0$& $           1175$ &     $  0$& $  0$& $  0$& $ -1$& $  2$\\
  $       -42581.735423$& $     -30.4$& $  -7$& $           1173$ &     $ -2$& $  4$& $  1$& $ -2$& $  0$\\
  $       -86549.627153$& $     -15.0$& $ -13$& $           1149$ &     $ -2$& $  3$& $  1$& $ -1$& $  0$\\
  $         1458.114127$& $     888.8$& $   0$& $           1103$ &     $ -2$& $  5$& $ -2$& $  0$& $  0$\\
  $        87963.797209$& $      14.7$& $  14$& $           1073$ &     $  0$& $  2$& $  0$& $ -1$& $  0$\\
  $       -22646.996366$& $     -57.2$& $  -3$& $           1060$ &     $  1$& $ -3$& $  0$& $  3$& $  0$\\
  $        21320.895355$& $      60.8$& $   3$& $           1031$ &     $  1$& $ -2$& $  0$& $  2$& $  0$\\
  $        87987.783148$& $      14.7$& $  14$& $            932$ &     $  0$& $  2$& $ -1$& $  0$& $  0$\\
  $      -173079.296175$& $      -7.5$& $ -27$& $            909$ &     $ -4$& $  6$& $  0$& $ -1$& $  0$\\
  $       305038.998861$& $       4.2$& $  47$& $            820$ &     $  4$& $ -3$& $  0$& $  0$& $  0$\\
  $       -63850.631098$& $     -20.3$& $ -10$& $            805$ &     $ -3$& $  6$& $  0$& $ -2$& $  0$\\
  $         2792.270745$& $     464.1$& $   0$& $            793$ &     $ -4$& $ 10$& $  0$& $ -5$& $  0$\\
  $       -19882.739378$& $     -65.2$& $  -3$& $            761$ &     $ -3$& $  7$& $  0$& $ -3$& $  0$\\
  $        65288.787089$& $      19.9$& $  10$& $            723$ &     $  1$& $ -1$& $  0$& $  1$& $  0$\\
  $        66670.915597$& $      19.4$& $  10$& $            711$ &     $ -1$& $  4$& $  0$& $ -2$& $  0$\\
  $          -56.106533$& $  -23098.9$& $   0$& $            679$ &     $  0$& $  0$& $ -1$& $  0$& $  2$\\
  $      -282307.961256$& $      -4.6$& $ -44$& $            655$ &     $ -5$& $  6$& $  0$& $  0$& $  0$\\
  $      -151762.428619$& $      -8.5$& $ -23$& $            599$ &     $ -3$& $  4$& $ -1$& $  1$& $  0$\\
  $      -151810.400498$& $      -8.5$& $ -23$& $            573$ &     $ -3$& $  4$& $  1$& $ -1$& $  0$\\
  $       174517.452169$& $       7.4$& $  27$& $            545$ &     $  2$& $ -1$& $  0$& $  0$& $  0$\\
  $       -43939.877972$& $     -29.5$& $  -7$& $            511$ &     $  0$& $ -1$& $  0$& $  2$& $  0$\\
  $       -65256.745538$& $     -19.9$& $ -10$& $            502$ &     $ -1$& $  1$& $  1$& $  0$& $  0$\\
  $      -107842.508766$& $     -12.0$& $ -17$& $            454$ &     $ -3$& $  5$& $  1$& $ -2$& $  0$\\
  $      -238340.069522$& $      -5.4$& $ -37$& $            450$ &     $ -5$& $  7$& $  0$& $ -1$& $  0$\\
  $         2768.284805$& $     468.2$& $   0$& $            430$ &     $ -4$& $ 10$& $  1$& $ -6$& $  0$\\
  $       370299.772209$& $       3.5$& $  57$& $            428$ &     $  5$& $ -4$& $  0$& $  0$& $  0$\\
  $       -21312.839748$& $     -60.8$& $  -3$& $            427$ &     $ -1$& $  2$& $  2$& $ -2$& $  0$\\
\trait
\end{tabular}
\caption{
Frequency analysis of $z_5=e_5\exp (i \varpi_5)$ in the Sun-Jupiter-Saturn system.
The integration of the complete Newton equations 
is performed over 50 myr with two output stepsizes  $h=200$ yr and $h'=200.0002$ yr.
$\nu_{0i}$ are the reconstructed frequencies using  formula
(\ref{eq.91},\ref{eq.92}) and the two step sizes $h$ and $h'$.
$P_i$ is the period of the terms, while $A_i$ their amplitude.
The different frequencies are identified 
as integer combination s of the fundamental frequencies $\nu_{0i} = k_{1i} n_5 + k_{2i} n_6 + k_{3i}  g_5 + k_{4i} g_6 + k_{5i} s_6$  }
\llabel{tab4}
\end{table}
\normalsize

\begin{table}
\small
\begin{tabular}{|  r| r|r|r| r|  r| r|r|r|}
\trait
$\nu_{0i}$ &            $P_i$ (years) &$k$& $\abs{A_i}\times 10^{8}$ &   $k_{1i}$&  $k_{2i}$&  $k_{3i}$&  $k_{4i}$&  $k_{5i}$ \\ 
\trait 
  $          -26.039362$& $  -49770.8$& $   0$& $         315418$ &     $  0$& $  0$& $  0$& $  0$& $  1$\\
  $            0.000000$& $          $& $   0$& $          59173$ &     $  0$& $  0$& $  0$& $  0$& $  0$\\
  $           82.066859$& $   15792.0$& $   0$& $           1457$ &     $  0$& $  0$& $  0$& $  2$& $ -1$\\
  $         1464.195357$& $     885.1$& $   0$& $            873$ &     $ -2$& $  5$& $  0$& $ -1$& $ -1$\\
  $           58.080919$& $   22313.7$& $   0$& $            820$ &     $  0$& $  0$& $  1$& $  1$& $ -1$\\
  $          -50.025302$& $  -25906.9$& $   0$& $            666$ &     $  0$& $  0$& $  1$& $ -1$& $  1$\\
  $           34.094979$& $   38011.5$& $   0$& $            623$ &     $  0$& $  0$& $  2$& $  0$& $ -1$\\
  $           -2.053422$& $ -631141.6$& $   0$& $            285$ &     $  0$& $  0$& $ -1$& $  1$& $  1$\\
  $        22757.076973$& $      56.9$& $   4$& $            228$ &     $ -1$& $  3$& $  0$& $  0$& $ -1$\\
  $          106.052799$& $   12220.3$& $   0$& $            197$ &     $  0$& $  0$& $ -1$& $  3$& $ -1$\\
  $         1440.209417$& $     899.9$& $   0$& $            168$ &     $ -2$& $  5$& $  1$& $ -2$& $ -1$\\
  $        65234.733984$& $      19.9$& $  10$& $             93$ &     $  1$& $ -1$& $  0$& $  0$& $  1$\\
  $        88017.850319$& $      14.7$& $  14$& $             92$ &     $  0$& $  2$& $  0$& $  0$& $ -1$\\
  $       -65286.812708$& $     -19.9$& $ -10$& $             82$ &     $ -1$& $  1$& $  0$& $  0$& $  1$\\
  $       -42503.696373$& $     -30.5$& $  -7$& $             74$ &     $ -2$& $  4$& $  0$& $  0$& $ -1$\\
  $        21266.842254$& $      60.9$& $   3$& $             71$ &     $  1$& $ -2$& $  0$& $  1$& $  1$\\
  $       130495.507329$& $       9.9$& $  20$& $             61$ &     $  2$& $ -2$& $  0$& $  0$& $  1$\\
  $       -21318.920978$& $     -60.8$& $  -3$& $             58$ &     $ -1$& $  2$& $  0$& $ -1$& $  1$\\
  $       153278.623665$& $       8.5$& $  24$& $             39$ &     $  1$& $  1$& $  0$& $  0$& $ -1$\\
  $        -1408.167857$& $    -920.3$& $   0$& $             38$ &     $  2$& $ -5$& $  0$& $  3$& $  1$\\
  $       218539.397009$& $       5.9$& $  34$& $             37$ &     $  2$& $  0$& $  0$& $  0$& $ -1$\\
  $        21242.856314$& $      61.0$& $   3$& $             34$ &     $  1$& $ -2$& $  1$& $  0$& $  1$\\
  $       -21234.800697$& $     -61.0$& $  -3$& $             32$ &     $ -1$& $  2$& $  1$& $  0$& $ -1$\\
  $          -74.011242$& $  -17510.9$& $   0$& $             30$ &     $  0$& $  0$& $  2$& $ -2$& $  1$\\
  $       -43993.931092$& $     -29.5$& $  -7$& $             28$ &     $  0$& $ -1$& $  0$& $  1$& $  1$\\
  $        43941.852367$& $      29.5$& $   7$& $             28$ &     $  0$& $  1$& $  0$& $ -1$& $  1$\\
  $         1512.167238$& $     857.0$& $   0$& $             27$ &     $ -2$& $  5$& $ -2$& $  1$& $ -1$\\
  $         1416.223477$& $     915.1$& $   0$& $             26$ &     $ -2$& $  5$& $  2$& $ -3$& $ -1$\\
  $       195756.280683$& $       6.6$& $  30$& $             26$ &     $  3$& $ -3$& $  0$& $  0$& $  1$\\
  $      -130547.586054$& $      -9.9$& $ -20$& $             25$ &     $ -2$& $  2$& $  0$& $  0$& $  1$\\
  $        66724.968703$& $      19.4$& $  10$& $             25$ &     $ -1$& $  4$& $  0$& $ -1$& $ -1$\\
  $        86527.615599$& $      15.0$& $  13$& $             24$ &     $  2$& $ -3$& $  0$& $  1$& $  1$\\
  $        22781.062913$& $      56.9$& $   4$& $             23$ &     $ -1$& $  3$& $ -1$& $  1$& $ -1$\\
  $         1356.089136$& $     955.7$& $   0$& $             19$ &     $ -2$& $  5$& $  0$& $ -3$& $  1$\\
  $           21.932518$& $   59090.3$& $   0$& $             19$ &     $  0$& $  0$& $ -2$& $  2$& $  1$\\
  $       283800.170356$& $       4.6$& $  44$& $             19$ &     $  3$& $ -1$& $  0$& $  0$& $ -1$\\
  $        44025.972649$& $      29.4$& $   7$& $             17$ &     $  0$& $  1$& $  1$& $  0$& $ -1$\\
  $      -107764.469717$& $     -12.0$& $ -17$& $             17$ &     $ -3$& $  5$& $  0$& $  0$& $ -1$\\
  $          130.038740$& $    9966.3$& $   0$& $             17$ &     $  0$& $  0$& $ -2$& $  4$& $ -1$\\
  $         2846.323859$& $     455.3$& $   0$& $             16$ &     $ -4$& $ 10$& $  0$& $ -4$& $ -1$\\
  $        22733.091033$& $      57.0$& $   4$& $             15$ &     $ -1$& $  3$& $  1$& $ -1$& $ -1$\\
  $       131985.742049$& $       9.8$& $  20$& $             15$ &     $  0$& $  3$& $  0$& $ -1$& $ -1$\\
  $       -63796.577987$& $     -20.3$& $ -10$& $             14$ &     $ -3$& $  6$& $  0$& $ -1$& $ -1$\\
  $       -19828.686253$& $     -65.4$& $  -3$& $             14$ &     $ -3$& $  7$& $  0$& $ -2$& $ -1$\\
  $       -22701.049474$& $     -57.1$& $  -4$& $             13$ &     $  1$& $ -3$& $  0$& $  2$& $  1$\\
  $       -21210.814756$& $     -61.1$& $  -3$& $             13$ &     $ -1$& $  2$& $  0$& $  1$& $ -1$\\
  $       -86579.694320$& $     -15.0$& $ -13$& $             12$ &     $ -2$& $  3$& $  0$& $ -1$& $  1$\\
  $        -1384.181921$& $    -936.3$& $   0$& $             12$ &     $  2$& $ -5$& $ -1$& $  4$& $  1$\\
  $       261017.054022$& $       5.0$& $  40$& $             12$ &     $  4$& $ -4$& $  0$& $  0$& $  1$\\
  $        44049.958588$& $      29.4$& $   7$& $             11$ &     $  0$& $  1$& $  0$& $  1$& $ -1$\\
\trait
\end{tabular}
\caption{ Same as Table \ref{tab4} for $\sin (i_5/2) \exp(i \Omega_5)$, 
where $i_5$ and $\Omega_5$ are the inclination and longitude of the node 
with respect to the ecliptic and equinox J2000 reference frame. 
The constant term is due to the invariance of the angular momentum.
It would not be present in the reference frame of the invariant plane, orthogonal to the 
angular momentum of the system.}
\llabel{tab5}
\end{table}

\subsection{Numerical examples}

\begin{table}
\begin{tabular}{| r|  r |}
\trait
     & (arcsec/years) \cr
\trait
$n_5$ & 109256.6788245339 \cr
$n_6$ &  43995.9054783976 \cr
$g_5$ &      4.0278083375 \cr
$g_6$ &     28.0137484932 \cr
$s_6$ &    -26.0393621745 \cr  
\trait
\end{tabular} 
\caption{Fundamental frequencies of the Sun-Jupiter-Saturn 
system, obtained by frequency analysis over 50 myr. }
\end{table} 

\subsubsection{Single periodic term}
In order to evaluate numerically the efficiency of the method described above, 
we first used a function $ f(t) = \exp(i \nu_0 \, t) $ with a single periodic 
term of frequency $\nu_0$. $f(t)$ is evaluated from $-1000$ to $+ 1000$, with a stepsize $h=1$,
and the frequency analysis is performed with different values of $\nu_0$, 
for values of $\nu_0/\pi$  higher than  $1000$. The value $\nu_0$ is only recovered 
for $\nu_0/\pi < 1$ (Table \ref{tab2}), but for larger values of 
$\nu_0$, the use of a second stepsize $h'=1.001$ and formula (\ref{eq.91},\ref{eq.92}) allows to recover  
the true frequency for values of $\nu_0/\pi$ as high as $1000$ with great accuracy  (Table \ref{tab3}). 

\subsubsection{Sun-Jupiter-Saturn system}

As a second, and more realistic example, we have consider a full numerical integration of the 
Sun-Jupiter-Saturn system over 50 millions of years. 
The numerical integration is performed with the symplectic integrator $SBAB_3$ that is well 
adapted to perturbed Hamiltonian systems (Laskar and Robutel, 2001).  In order to illustrate the 
previous section, we have performed the frequency analysis 
of the variable $z_5=e_5\exp (i \varpi_5)$ (Fig.\ref{tab4}), where $e_5$ is the eccentricity of Jupiter, and 
$\varpi_5$ its longitude of perihelion, and the analysis of 
$\sin (i_5/2) \exp(i \Omega_5)$ (Fig.\ref{tab5}), 
where $i_5$ and $\Omega_5$ are the inclination and longitude of the node 
with respect to the ecliptic and equinox J2000 reference frame. 
If one considers the invariance of the angular momentum, this problems 
has  5 degrees of freedom, that should correspond to 5  fundamental frequencies for a regular 
KAM solution. These frequencies will be $n_5, n_6$, the mean mean motion of Jupiter and Saturn,
$g_5, g_6$ related to the precessional motion of the perihelion of Jupiter and Saturn, 
and $s_6$ related to the motion of the node of the two orbits (see Laskar, 1990).
In the integration, we used an output  stepsize of $200$ years, while the integration stepsize
is $0.1$ year. This allows to recover precisely the secular frequencies of the system, 
but as the short period perturbations are important,
many terms appear in the quasiperiodic decomposition that are in fact aliased 
terms coming from the short period terms (see Table 4,5).

When we use an additional  output time   $h'=200.0002$ years, and formula (\ref{eq.91},\ref{eq.92}),
we can recover  the true value of the aliased frequencies $\nu_{0i}$, even if 
some of the periods are smaller than 6 years (columns $P_i$). In Table 4, 5, 
$k$ is the integer appearing in formula (\ref{eq.84}), denoting the number 
of turns that have been "lost" by  the large stepsize.
All the determined frequencies $\nu_{0i}$ can then be identified as integer combinations 
\be
\nu_{0i} = k_{1i} n_5 + k_{2i} n_6 + k_{3i}  g_5 + k_{4i} g_6 + k_{5i} s_6
\ee 
of the fundamental frequencies $n_5, n_6, g_5, g_6, s_6$ given  in Table 5.
The full solution can then be compared to the quasiperiodic solution 
obtained by iteration of (Bretagnon and Simon, 1990).

The presence of a constant term in Table \ref{tab5} is due to the invariance of the angular momentum.
It would not be present in the reference frame of the invariant plane, orthogonal to the 
angular momentum of the system (see Malige \etal, 2002).

The existence of short periods is thus not an obstacle, and the aliasing problem 
can be overcome, as in practice, using two 
output time steps is not a big constraint. It should not be costly in term 
of CPU time, as the two output are generated during the same integration. 
The large advantage of this procedure versus the numerical averaging 
that was usually performed online (as in Nobili \etal, 1989), 
is that no information is lost, and the whole solution can be recovered.

\thispagestyle{empty}

\vspace{0.5 cm}
{\bf Acknowledgements} The author is very obliged to 
A. Chenciner, F. Joutel, L. Niederman and D. Sauzin for very useful discussions
 at several stages of this work. 
 The help of F. Joutel and P. Robutel for the numerical comparisons 
is deeply appreciated.

\end{document}